\documentclass[12pt,a4paper]{article}

\usepackage{amsfonts}
\usepackage{amsmath, amsthm}
\usepackage{latexsym}
\usepackage{amssymb}
\usepackage{amscd}
\usepackage{verbatim}
\usepackage{cite}

\makeatletter
\@addtoreset{equation}{section}
\makeatother

\newtheorem{theorem}{Theorem}[section]

\theoremstyle{definition}
\newtheorem{definition}[theorem]{Definition}
\newtheorem{example}[theorem]{Example}
\newtheorem{remark}[theorem]{Remark}

\begin{document}

\title{\textbf{Tables of parameters of symmetric configurations }$v_{k}$\footnote{This
research was supported by the Italian Ministry MIUR,
 Geometrie di Galois e strutture di incidenza, PRIN
2009--2010 and by the INdAM group G.N.S.A.G.A.}}
\author{Alexander A. Davydov \\
Institute for Information Transmission Problems (Kharkevich institute),\\
Russian Academy of Sciences, Bol'shoi Karetnyi per. 19, GSP-4, \\
Moscow, 127994, Russian Federation, e-mail: adav@iitp.ru \and Giorgio Faina,
Massimo Giulietti, \and Stefano Marcugini and Fernanda Pambianco \\
Dipartimento di Matematica e Informatica, \\
Universit\`{a} degli Studi di Perugia, Via Vanvitelli 1,
Perugia, 06123,
Italy\\
e-mail: \{faina, giuliet, gino, fernanda\} @dmi.unipg.it}
\date{}
\maketitle

\begin{abstract}
Tables of the currently known parameters of symmetric
configurations are given. Formulas for parameters of the known
infinite families of symmetric configurations are presented as
well. The results of the recent paper \cite{DFGMP-submitted}
are used. This work can be viewed as
 an appendix to \cite{DFGMP-submitted}, in the sense that the tables
 given here cover a much larger set of parameters.
\end{abstract}



\section{Introduction\label{intro}}

Configurations as combinatorial structures were defined in
1876. For an introduction to the problems connected with
 configurations, see \cite{Gropp-ConfGeomCombin,Gropp-Handb,Grunbaum} and the
references therein.

\begin{definition}
\label{Def1-Configur}\cite{Gropp-Handb}

\textbf{(i)} A configuration $(v_{r},b_{k})$ is an incidence structure of $v$
points and $b$ lines such that each line contains $k$ points, each point
lies on $r$ lines, and two distinct points are connected by \emph{at most}
one line.

\textbf{(ii)} If $v=b$ and, hence, $r=k$, the configuration is \emph{\
symmetric}, and it is referred to as a configuration $v_{k}$.

\textbf{(iii)} The \emph{deficiency} $d$ of a configuration $(v_{r},b_{k})$
is the value $d=v-r(k-1)-1$.
\end{definition}

A symmetric configuration $v_k$ is \emph{cyclic} if there exists a
permutation of the set of its points mappings blocks to blocks, and acting
regularly on both points and blocks. Equivalently, $v_k$ is cyclic if one of
its incidence matrix is circulant.

Steiner systems are configurations with $d=0$
\cite{Gropp-Handb}. The \emph{deficiency} of a symmetric
configuration $v_{k}$ is $d=v-(k^{2}-k+1)$. The deficiency of
$v_{k}$ is zero if and only if $v_{k}$ is a finite projective
plane of order $k-1$.

A configuration $(v_{r},b_{k})$ can be viewed as a $k$-uniform
$r$ -regular linear hypergraph with $v$ vertices and $b$
hyperedges \cite{Gropp-ConfGraph,Gropp-Handb}. Connections of
configurations $(v_{r},b_{k})$ with numerical semigroups are
noted in \cite{BrasAmorosSemigroup,StokesBrasAmoros2013}. Some
analogies between configurations $(v_{r},b_{k})$, regular
graphs, and molecule models of chemical elements are remarked
in \cite{Gropp-Chemic}. As an example of a practical
application of configurations (both symmetric and
non-symmetric), we mention also the problem of user privacy for
using database; see
\cite{DomingoAmorosPeertoPeer,StokesPeertoPeer} and the
references therein.

Denote by $\mathbf{M}(v,k)$ an incidence matrix of a symmetric configuration
$v_{k}.$ Any matrix $\mathbf{M}(v,k)$ is a $v\times v$ 01-matrix with $k$
ones in every row and column; moreover, the $2\times 2$ matrix consisting of
all ones is not a submatrix of $\mathbf{M}(v,k)$. Two incidence matrices of
the same configuration may differ by a permutation on the rows and the
columns.

A matrix $\mathbf{M}(v,k)$ can be considered as a biadjacency
matrix of the \emph{Levi graph} of the configuration $v_{k}$
which is a $k$-regular bipartite graph without multiple edges
\cite[Sec.\thinspace 7.2]{Gropp-Handb}. Clearly, the graph has
girth at least~six, i.e. it does not contain 4-cycles. Such
graphs are useful for the construction of bipartite-graph codes
that can be treated as\emph{\ low-density parity-check}
(LDPC)\emph{\ codes}. If $\mathbf{M}(v,k)$ is \emph{circulant},
then the corresponding LDPC code is \emph{quasi-cyclic}; it can
be encoded with the help of shift-registers with relatively
small complexity; see
\cite{AfDaZ,AfDaZ-InfProc,DGMP-ACCT2008,DGMP-GraphCodes,GabidISIT,%
CyclDecomp,CyclDecomp2,QCEncoder} and the references therein.

Matrices $\mathbf{M}(v,k)$ consisting of square circulant
submatrices have a number of useful properties, e.g. they are
more suitable for LDPC codes implementation. We say that a
01-matrix\emph{\ }$\mathbf{A}$ is \emph{block double-circulant}
(BDC for short) if $\mathbf{A }$ consists of square circulant
blocks whose weights give rise to a circulant matrix (see
Definition~\ref{def3.1_block double-circul}). A configuration
$v_{k}$ with a BDC incidence matrix $\mathbf{M} (v,k) $ is
called a \emph{BDC symmetric configuration.} Symmetric and
non-symmetric configurations with incidence matrices consisting
of square circulant blocks are considered, for example, in
\cite{AFLN-Semiplanes,AFLN-CyclSchem,DGMP-ACCT2008,DGMP-Petersb2009,DGMP-GraphCodes,CyclDecomp,CyclDecomp2,Pepe}.
In \cite{AFLN-Semiplanes,AFLN-CyclSchem}, BDC symmetric
configuration are considered in connection with $\mathbb
Z_\mu$-schemes (see Remarks \ref{zetamu} and \ref{zetamu2} in
Section \ref{sec_bdc&families}).

Cyclic configurations are considered, for instance, in
\cite{AFLN-Semiplanes,AFLN-CyclSchem,DGMP-ACCT2008,DGMP-Petersb2009,DGMP-GraphCodes,Funk2008,Gropp-nk,CyclDecomp,CyclDecomp2,Lipman,MePaWolk}.
A standard method to construct cyclic configurations (or,
equivalently, circulant matrices $\mathbf{M}(v,k)$) is based on
\emph{Golomb rulers}
\cite{Dimit,Funk2008,Gropp-nk,Shearer-Handb,ShearerWebShortest,ShearerWebModulGR}.

\begin{definition}
\label{Def1_GR}\cite{Shearer-Handb,Funk2008}

\textbf{(i)} A \emph{Golomb ruler} $\mathrm{G}_{k}$ of
\emph{order} $k$ is an ordered set of $k$ integers
$(a_{1},a_{2},\ldots ,a_{k})$ such that $ 0\leq
a_{1}<a_{2}<\ldots <a_{k}$ and all the differences $
\{a_{i}-a_{j}\,|\,1\leq j<i\leq k\}$ are distinct. The
\emph{length} $L_{ \mathrm{G}}(k)$ of the ruler
$\mathrm{G}_{k}$ is equal to $a_{k}-a_{1}$. Let
$L_{\overline{\mathrm{G}}}(k)$ be the length of the
\emph{shortest known} Golomb ruler $\overline{\mathrm{G}}_{k}$.

\textbf{(ii)} A $(v,k)$ \emph{modular Golomb ruler} is an ordered set of $k$
integers $(a_{1},a_{2},\ldots ,a_{k})$ such that $0\leq a_{1}<a_{2}<\ldots
<a_{k}$ and all the differences $\{a_{i}-a_{j}\,|\,1\leq i,j\leq k;$ $i\neq
j\}$ are distinct and nonzero modulo $v$.
\end{definition}

For any $\delta \geq 0$, Golomb rulers $(a_{1},a_{2},\ldots
,a_{k})$ and $ (a_{1}+\delta ,a_{2}+\delta ,\ldots
,a_{k}+\delta )$ have the same properties. Usually, $a_{1}=0$
is assumed. We say that a 0,1-\emph{vector}
${\mathbf{u}}=(u_{0},u_{1,}\ldots ,u_{v-1})$ \emph{corresponds}
to a (modular) Golomb ruler if the increasing sequence of
integers $j\in \{0,1,\ldots,v-1\}$ such that $u_j=1$ form a
(modular) Golomb ruler.

Recall that \emph{weight} of a \emph{circulant }$0,1$\emph{-matrix }is the
number of ones in each its row.

\begin{theorem}
\label{Th1_2L+1}\emph{\cite[Sec.\thinspace 4]{Gropp-nk},\cite{Longyear1975}}

\emph{\textbf{(i)}} Any Golomb ruler $\mathrm{G}_{k}$ of length
$L_{\mathrm{ G }}(k)$ is a $(v,k)$ modular Golomb ruler for all
$v$ such that $v\geq 2L_{ \mathrm{G} }(k)+1$.

\emph{\textbf{(ii)}} A circulant $v\times v$ \emph{0,1}-matrix of weight $k$
is an incidence matrix $\mathbf{M}(v,k)$ of a cyclic symmetric configuration
$v_{k} $ if and only if the first row of the matrix corresponds to a $(v,k)$
modular Golomb ruler.

\emph{\textbf{(iii)}} For all \/\/$v$ such that $v\geq
2L_{\overline{\mathrm{ G}}}(k)+1,$ there exists a cyclic
symmetric configuration$~v_{k}$.
\end{theorem}

We call the value $G(k)=2L_{\overline{\mathrm{G}}}(k)+1$ the
\emph{Golomb bound}. On the other hand, we call $
P(k)=k^{2}-k+1$ the \emph{projective plane bound}.

Let $v_{\delta }(k)$ be the smallest possible value of $v$ for
which a $ (v,k) $ modular Golomb ruler (or, equivalently, a
cyclic symmetric configuration) exists.

In \cite{DFGMP-submitted}, two bounds are considered. The
\emph{existence bound} $E(k)$ is the least
 integer such that for any $v\geq E(k)$, there exists a
symmetric configuration $v_{k}$. Similarly, the \emph{cyclic
existence bound} $E_{c}(k)$ is the least integer such that for
any $v\geq E_{c}(k)$, there exists a cyclic~$v_{k}$. Clearly,
for a fixed $k$, we have
\begin{eqnarray}
k^{2}-k+1 &=&P(k)\leq E(k)\leq E_{c}(k)\leq G(k)=2L_{\overline{\mathrm{G}}
}(k)+1.  \label{eq1_bounds} \\
k^{2}-k+1 &=&P(k)\leq v_{\delta }(k)\leq E_{c}(k)\leq G(k)=2L_{\overline{
\mathrm{G}}}(k)+1.  \label{eq1_boundsCyclic}
\end{eqnarray}

The aim of this work is to give tables of the currently known
parameters of symmetric configurations $v_{k},$ including those
arising from the
 recent work \cite{DFGMP-submitted}. We consider the\emph{\
spectrum }of possible {parameters} of $v_{k}$ (with special
attention to {\ parameters} of {cyclic symmetric
configurations}) in the interval
\begin{equation}
k^{2}-k+1=P(k)\leq v<G(k)=2L_{\overline{\mathrm{G}}}(k)+1.
\label{eq1_region}
\end{equation}
Also, we pay attention to {parameters} of {circulant} and block
double-circulant incidence matrices $\mathbf{M}(v,k)$. Some
upper bounds on $ E(k)$ and $E_{c}(k)$ are pointed out.

From the stand point of applications, including Coding Theory,
it is sometimes useful to have different matrices
${\mathbf{M}}(v,k)$ for the same $v$ and $k$. This is why we
remark situations when \emph{different constructions} provide
configurations with \emph{the same parameters}.

The \emph{Generalized Martinetti Construction} (Construction
GM) proposed in \cite{Funk1993}
  plays a key role for the
investigation of the spectrum of possible parameters of
symmetric configurations as it provides, for a fixed $k$,
intervals of values of $v$ for which a $v_k$ exists.
Construction GM has been considered also in
\cite{AFLN-Semiplanes,AfDaZ,AfDaZ-InfProc}\footnote{ The
authors of the papers \cite{AfDaZ,AfDaZ-InfProc} (represented
here by Davydov) regret that  the paper \cite{Funk1993} is not
cited in \cite{AfDaZ,AfDaZ-InfProc}; the reason is that,
unfortunately, the authors did not know the paper
\cite{Funk1993} during the preparation of
\cite{AfDaZ,AfDaZ-InfProc}.} To be successfully applied,
 Construction GM needs a convenient starting incidence
matrix. To this end, BDC matrices turn out to be particularly
useful. In this work  new starting matrices proposed in
\cite{DFGMP-submitted} are considered as well as those
originally proposed in \cite{AFLN-Semiplanes,AfDaZ-InfProc}.

We remark that new cyclic configurations provide new modular Golomb rulers,
i.e. new deficient cyclic difference~sets.

The work is organized as follows. In Section \ref{sec_known},
we briefly summarize some constructions and parameters of
configurations $v_{k}$. Preliminaries on BDC matrices are given
in Section \ref{sec_bdc&families}. In Section \ref{sec_BDC},
parameters of block double-circulant incidence matrices
$\mathbf{M}(v,k)$ are reported, according to some results from
\cite{DFGMP-submitted}. In Section~\ref{sec_admitExten},
parameters of configurations $v_{k}$ obtained by
 the Construction
GM from the starting matrices proposed in
\cite{AFLN-Semiplanes,AfDaZ-InfProc,DFGMP-submitted} are given.
In Sections \ref{sec_ParamCyclicConfig}
and~\ref{sec6_spectrum}, results on the spectra of parameters
of cyclic and non-cyclic configurations are reported. Finally,
Section \ref{sec_tables} contains the tables of parameters of
symmetric configurations which are the main object of the
paper. In particular, tables of parameters of\ BDC\
configurations $v_{k}$ based on projective planes and punctured
 affine planes are given,
 as well as tables of values $v$ for
which a cyclic symmetric configuration $v_{k}$ exists. Finally,
aggregated tables on the existence of symmetric configurations
are given. In the tables, the new parameters obtained from
\cite{DFGMP-submitted} are written in bold font.

\section{Some known constructions and parameters of configurations $v_{k}$
with $P(k)\leq v<G(k)$\label{sec_known}}

The aim of this section is to provide a list of pairs $(v,k)$
for which a (cyclic) symmetric configuration $v_{k}$ is known
to exist, see Equations
\eqref{eq2_cyclicPG(2,q)}--(\ref{eq2_AaParBalb}). Infinite
families of configurations $v_{k}$ given in this section are
considered in
\cite{AFLN-graphs,AFLN-ConfigGraphs,AFLN-Semiplanes,AFLN-CyclSchem,AfDaZ,%
AfDaZ-InfProc,ArParBalbNetw2011,ArParBalbHegDM2010,BalbuenaSIAM2008,Bose,DFGMPConfigArxiv2012,%
DGMP-ACCT2008,DGMP-Petersb2009,DGMP-GraphCodes,Funk1993,Funk2008,FunkLabNap,GH,Gropp-nk,%
Gropp-Chemic,Gropp-nonsim,Gropp-ConfGraph,Gropp-ConfGeomCombin,Gropp-Handb,Grunbaum,Lipman,%
MePaWolk,Ruzsa,Shearer-Handb,Singer}; see also the references
therein.

Throughout the work, $q$ is a prime power and $p$ is a prime. Let $F_{q}$ be
Galois field of $q$ elements. Let $F_{q}^{\ast }=F_{q}\backslash \{0\}.$ Let
$\mathbf{0}_{u}$ be the zero $u\times u$ matrix. Denote by $\mathbf{P}_{u}$
a permutation matrix of order $u.$

We recall that several pairs $(v,k-\delta)$ can be actually
obtained from a given $v_k$; it is a basic result on symmetric
configurations.

\begin{theorem}
\label{Th2-Mend}
\emph{\cite[Sec.\thinspace2]{AFLN-Semiplanes},\cite[Sec.\thinspace5.2]{Gropp-ConfGraph},\cite[Sec.~2.5]{Grunbaum}\cite{MePaWolk}}
If a (cyclic) configuration $v_{k}$ exists, then for each
$\delta $ with $0\leq \delta <k$ there exists a (cyclic)
configuration $v_{k-\delta }$ as well.
\end{theorem}

We note that from
 a cyclic configuration $v_{k}$, a
cyclic configurations $v_{k-\delta }$ can be obtained by
dismissing $\delta $ ones in the 1-st row of its incidence
matrix. For the general case, Theorem~\ref{Th2-Mend} is based
on the fact that an incidence matrix $\mathbf{M}(v,k)$ can be
represented as a sum of $k$ permutations $v\times v$ matrices
(in different ways). This fact follows from the results of
Steinits (1894) and K\"{o}nig (1914), see e.g.
\cite[Sec.\thinspace 5.2]{Gropp-ConfGraph} and
\cite[Sec.~2.5]{Grunbaum}.

The value $\delta $ appearing in Equations
\eqref{eq2_cyclicPG(2,q)}--(\ref{eq2_AaParBalb}) is connected
with Theorem \ref{Th2-Mend}. When a reference is given, it
usually refers to the case $\delta =0$.

The families giving rise to pairs
(\ref{eq2_cyclicPG(2,q)})--(\ref{eq2_cyclicRuzsa}) below are
obtained from $(v,k)$ modular Golomb rulers \cite[Ch.\thinspace
5]{Dimit}, \cite{Draka}, \cite[Sec.\thinspace 5]{Gropp-nk},
\cite[Sec.\thinspace 19.3]{Shearer-Handb}; see Theorem
\ref{Th1_2L+1}(ii).
\begin{eqnarray}
\text{cyclic }v_{k} &:&v=q^{2}+q+1,\hspace{0.2cm}k=q+1-\delta ,\hspace{0.2cm}
q+1>\delta \geq 0;  \label{eq2_cyclicPG(2,q)} \\
\text{cyclic }v_{k} &:&v=q^{2}-1,\hspace{0.2cm}k=q-\delta ,\hspace{0.2cm}
q>\delta \geq 0;  \label{eq2_cyclicAG(2,q)} \\
\text{cyclic }v_{k} &:&v=p^{2}-p,\hspace{0.2cm}k=p-1-\delta ,\hspace{0.2cm}
p-1>\delta \geq 0.  \label{eq2_cyclicRuzsa}
\end{eqnarray}
The configurations giving rise to (\ref{eq2_cyclicPG(2,q)}) use
the incidence matrix of the cyclic projective plane $PG(2,q)$
\cite[Sec.\thinspace 5.5]{Dimit}, \cite{Draka},
\cite[Th.\thinspace 19.15]{Shearer-Handb},\cite{Singer}. The
family with parameters (\ref{eq2_cyclicAG(2,q)}) is obtained
from the \emph{cyclic punctured affine plane} $AG(2,q)$
\cite{Bose}, \cite[Sec.\thinspace 5.6]{Dimit},
\cite{Draka,Funk2008}, \cite[Th.\thinspace
19.17]{Shearer-Handb}; see also \cite[Ex.\thinspace
5]{DGMP-GraphCodes} and \cite{FunkLabNap} where the
configurations are called \emph{anti-flags}. We recall that the
punctured plane $AG(2,q)$ is the affine plane without the
origin and the lines through the origin. Punctured affine
planes are also called elliptic (Desarguesian) semiplanes of
type L. Finally, the configurations with parameters
(\ref{eq2_cyclicRuzsa}) follow from Ruzsa's construction
\cite[Sec.\thinspace 5.4]{Dimit}, \cite{Draka,Ruzsa},
\cite[Th.\thinspace 19.19]{Shearer-Handb}.

The non-cyclic families with parameters (\ref{eq2_q-cancel})
and (\ref{eq2_q-1-cancel}) are given in \cite[Constructions
(i),(ii), p.\thinspace 126 ]{AFLN-graphs} and
\cite[Constructions 3.2,3.3, Rem.\thinspace 3.5]{GH}; see also
the references therein and \cite{AfDaZ}, \cite[Sec.\thinspace
3]{AfDaZ-InfProc}, \cite[Sec.\thinspace 7.3]{DGMP-GraphCodes},
\cite{FunkLabNap,Gropp-nk,MePaWolk}.
\begin{eqnarray}
v_{k} &:&v=q^{2}-qs,\hspace{0.2cm}k=q-s-\delta ,\hspace{0.2cm}q>s\geq 0,
\hspace{0.2cm}q-s>\delta \geq 0;  \label{eq2_q-cancel} \\
v_{k} &:&v=q^{2}-(q-1)s-1,\hspace{0.2cm}k=q-s-\delta ,\hspace{0.2cm}q>s\geq
0,\hspace{0.2cm}q-s>\delta \geq 0.  \label{eq2_q-1-cancel}
\end{eqnarray}

For $q$ a \emph{square}, in \cite[Conjec.\thinspace 4.4,
Rem.\thinspace 4.5, Ex.\thinspace 4.6]{AFLN-graphs},
\cite[Th.\thinspace 6.4]{AFLN-Semiplanes}, and
\cite[Construction 3.7, Th.\thinspace 3.8]{GH}, families of
non-cyclic configuration $v_{k}$ with parameters
(\ref{eq2_Baer}) are provided; see also \cite[Ex.\thinspace
8]{DGMP-GraphCodes}. Taking $c=q- \sqrt{q}$, we obtain
(\ref{eq2_Baer2}), see also \cite[ Ex.\thinspace
2(ii)]{DGMP-ACCT2008}, \cite{FunkLabNap}.
\begin{eqnarray}
v_{k} &:&v=c(q+\sqrt{q}+1),k=\sqrt{q}+c-\delta ,\,c=2,3,\ldots,q-\sqrt{q}
,\,\delta \geq 0;  \label{eq2_Baer} \\
v_{k} &:&v=q^{2}-\sqrt{q},\hspace{0.2cm}k=q-\delta ,\hspace{0.2cm}q>\delta
\geq 0.  \label{eq2_Baer2}
\end{eqnarray}

In \cite[Th.\thinspace 1.1]{FunkLabNap}, a non-cyclic family with parameters
\begin{equation}
v_{k}:v=2p^{2},\hspace{0.2cm}k=p+s-\delta ,\hspace{0.2cm}0<s\leq q+1,\hspace{
0.2cm}\,q^{2}+q+1\leq p,\hspace{0.2cm}p+s>\delta \geq 0
\label{eq2_FuLabNabDecomp}
\end{equation}
is given. In \cite[Sec.\thinspace \thinspace 6]{DGMP-GraphCodes}, a
construction of non-cyclic configuration based on the cyclic punctured affine
plane, is provided with parameters
\begin{eqnarray}
v_{k} &:&v=c(q-1),\hspace{0.2cm}k=c-\delta ,\hspace{0.2cm}c=2,3,\ldots ,b,
\hspace{0.2cm}b=q\,\,\mathrm{if}\,\,\delta \geq 1,\,\,
\label{eq2_affine_q-1} \\
b &=&\left\lceil \frac{q}{2}\right\rceil \,\,\mathrm{if}\,\,\delta =0,\text{
}c>\delta \geq 0.  \notag
\end{eqnarray}

In \cite[Sec.\thinspace 2]{DGMP-ACCT2008}, \cite[Sec.\thinspace
\thinspace 3]{DGMP-GraphCodes}, the following geometrical
construction which uses point orbits under the action of a
collineation group is described.\smallskip

\noindent\textbf{Construction A}. Take any point orbit
$\mathcal{P}$ under the action of a collineation group in an
affine or projective space of order $q$. Choose an integer
$k\leq q+1$ such that the set $\mathcal{L}(\mathcal{P} ,k)$ of
lines meeting $\mathcal{P}$ in precisely $k$ points is not
empty. Define the following incidence structure: the points are
the points of $ \mathcal{P}$, the lines are the lines of
$\mathcal{L}(\mathcal{P},k)$, the incidence is that of the
ambient space.

\begin{theorem}
\label{Th3_constrOneOrb} In Construction A the number of lines
of $\mathcal{ \ L }(\mathcal{P},k)$ through a point of
$\mathcal{P}$ is a constant $r_{k}$. The incidence structure is
a configuration $(v_{r_{k}},b_{k})$ with $
v_{r_{k}}=|\mathcal{P}|$, $b_{k}=|\mathcal{L}(\mathcal{P},k)|$.
\end{theorem}

By Definition \ref{Def1-Configur}, if $r_{k}=k$ then Construction A
produces a symmetric configuration $v_{k}$. It is noted in
\cite{DGMP-ACCT2008,DGMP-GraphCodes} that Construction A works
for any $2$-$ (v,k,1)$ design $D$ and for any group of
automorphism of $D$. The size of any block in $D$ plays the
role of $q+1$.

Families of non-cyclic configuration $v_{k}$ obtained by
Construction A with the following parameters are given in
\cite[Exs\thinspace 2,\thinspace 3]{DGMP-GraphCodes}.
\begin{eqnarray}
v_{k} &:&v=\frac{q(q-1)}{2},\hspace{0.2cm}k=\frac{q+1}{2}-\delta ,\hspace{
0.2cm}\frac{q+1}{2}>\delta \geq 0,\text{ }q\text{ odd}.
\label{eq2_ConstrA_GraphCode1} \\
v_{k} &:&v=\frac{q(q+1)}{2},\hspace{0.2cm}k=\frac{q-1}{2}-\delta ,\hspace{
0.2cm}\frac{q-1}{2}>\delta \geq 0,\text{ }q\text{ odd}.
\label{eq2_ConstrA_GraphCode2} \\
v_{k} &:&v=q^{2}+q-q\sqrt{q},\text{ }k=q-\sqrt{q},\text{ }q-\sqrt{q}>\delta
\geq 0,\text{ }q\text{ square.}  \label{eq2_ConstrA_GraphCode3}
\end{eqnarray}

In
\cite{ArParBalbNetw2011,ArParBalbHegDM2010,BalbuenaSIAM2008},
non-cyclic families with parameters
\eqref{eq2_Balbuena1}--\eqref{eq2_AaParBalb} are described  in
connection with graph theory; see also \cite{AFLN-CyclSchem}
for another construction of \eqref{eq2_Balbuena2}.
\begin{eqnarray}
v_{k} &:&v=q^{2}-rq-1,\hspace{0.2cm}k=q-r-\delta,\hspace{
0.2cm}q-r>\delta \geq 0,\text{ }q-3\ge r\ge0.
\label{eq2_Balbuena1} \\
v_{k} &:&v=q^{2}-q-2,\hspace{0.3cm}k=q-1-\delta ,\hspace{
0.2cm}q-1>\delta \geq 0.\label{eq2_Balbuena2} \\
v_{k} &:&v=tq-1,\hspace{0.87cm}k=t-\delta,\hspace{0.87cm}t>\delta \geq 0,\quad q> t\ge3.  \label{eq2_AaParBalbHeg}\\
v_{k} &:&v=tq-2,\hspace{0.87cm}k=t-\delta,\hspace{0.87cm}t>\delta \geq 0,\quad q> t\ge3.  \label{eq2_AaParBalb}
\end{eqnarray}

A classical construction by V. Martinetti for configurations
$v_{3}$, going back to 1887 \cite{Martinetti}, is described in
detail, e.g. in
\cite{AFLN-Semiplanes,AfDaZ-InfProc,Boben-v3,CarsDinStef-Reduc-n3,Funk1993,Gropp-Chemic},
\cite[ Sec. 2.4, Fig. 2.4.1]{Grunbaum}. In \cite{Funk1993} a
\emph{Generalized Martinetti Construction} (Construction GM)
for configurations $v_{k}$, $ k\geq 3$, is proposed. Use of
Construction GM to obtain a wide spectrum of symmetric
configurations parameters is considered in
\cite{AFLN-Semiplanes,AfDaZ,AfDaZ-InfProc}.

Construction GM can be presented from different points of view,
see \cite{AFLN-Semiplanes,Funk1993}. Here we focus on an
approach based on incidence matrices, which will be used for
obtaining new values of $v,k$.

\begin{definition}
\label{def2_ExtenAgeg} Let $\mathbf{\ M}(v,k)$ be an incidence
matrix of a symmetric configuration $v_{k}.$ In
$\mathbf{M}(v,k)$, we consider an aggregate $\mathcal{A}$ of
$k-1$ rows corresponding to pairwise disjoint lines of $v_{k}$
and $k-1$ columns corresponding to pairwise non-collinear
points of $v_{k}$. If a $(k-1)\times (k-1)$ submatrix
$\mathbf{C}{(\mathcal{A })}$ formed by the intersection of the
rows and columns of $\mathcal{A}$ is a permutation matrix
$\mathbf{P}_{k-1}$ then $\mathcal{A}$ is called an
\emph{extending aggregate} (or \emph{E-aggregate}). The matrix
$\mathbf{M} (v,k)$ \emph{admits an extension} if it contains at
least one E-aggregate. The matrix $\mathbf{M}(v,k)$
\emph{admits }$\theta $\emph{\ extensions} if it contains
$\theta $ E-aggregates that do not intersect each other.
\end{definition}

\noindent \textbf{Procedure E }(\emph{Extension Procedure}). Let $\mathbf{M}
(v,k)=[m_{ij}]$ be an incidence matrix of a symmetric configuration $v_{k}=(
\mathcal{P},\mathcal{L}).$ Assume that $\mathbf{M}(v,k)$ admits an
extension. A matrix $\mathbf{M}(v+1,k)$ can be obtained by two steps.

1. To the matrix $\mathbf{M}(v,k)$, add a new row from below and a new
column to the right. Denote the new $(v+1)\times (v+1)$ matrix by $\mathbf{B}
=[b_{ij}]$, and let $b_{v+1,v+1}=1$ whereas $b_{v+1,1}=\ldots =b_{v+1,v}=0,$
$b_{1,v+1}=\ldots =b_{v,v+1}=0.$

2. One of E-aggregates of $\mathbf{M}(v,k)$, say $\mathcal{A},$
is chosen. In the matrix $\mathbf{B,}$ we \textquotedblleft
clone\textquotedblright\ all $k-1$ ones of the submatrix
$\mathbf{C }{(\mathcal{A})}$ writing their \textquotedblleft
projections\textquotedblright\ to the new row and column.
Finally, the ones cloned are changed by zeroes. In other words,
let the aggregate $\mathcal{A}$ consist of rows with indexes
$i_{u},$ $u=1,2,\ldots ,k-1,$ and columns with indexes $j_{d},$
$d=1,2,\ldots ,k-1.$ Then the ones of
$\mathbf{C}{(\mathcal{A})}$ are as follows:
$m_{i_{u}j_{\pi(u)}}=1,$ $ u=1,2,\ldots ,k-1,$ for some
permutation $\pi$ of the indexes $1,\ldots,k-1$ . Then
$\mathbf{B}$ arising from Step 1 is changed as follows: $
b_{i_{u},v+1}=1,$ $b_{v+1,j_{d}}=1,$ $b_{i_{u}j_{\pi(u)}}=0,$
$u=1,2,\ldots ,k-1$, $d=1,2,\ldots ,k-1$.

As the Golomb bound $2L_{\overline{\mathrm{G}}}(k)+1$ is
important for studying parameters $v,k$, we note that for
sufficiently large orders $k$, relatively short Golomb rulers
are constructed and are available online, see
\cite{Dimit,ShearerWebShortest,ShearerWebModulGR} and the
references therein. For $k\leq 150$, the order of magnitude of
the lengths $L_{ \overline{\mathrm{G}}}(k)$ of the shortest
known Golomb rulers is $ck^{2}$ with $c\in \lbrack 0.7,0.9]$,
see
\cite{Dimit,Funk2008,Gropp-nk,Gropp-Handb,Shearer-Handb,ShearerWebShortest}.
Moreover, $L_{\overline{\mathrm{G}}}(k)<k^{2}\text{ for
}k<65000,$ see \cite{Dimit}. Constructions of Golomb rulers for
large $k$ can be found in \cite{Draka}. Remind also that Sidon
sets are equivalent to Golomb rulers, see
\cite{Dimit,BibliogrSidon} and the references therein.

\section{Preliminaries on block double-circulant incidence matrices \label{sec_bdc&families}}

Results of this section are taken from \cite{DFGMP-submitted},
see also the references therein, in particular, discussions of
\cite[Rem. 1,2]{DFGMP-submitted}.

Recall that the weight of a circulant binary is the number of $1$'s in each
its rows.

\begin{definition}
\label{def3.1_block double-circul}Let $v=td.$ A binary $v\times
v$ matrix $ \mathbf{A}$ is said to be a \emph{block
double-circulant matrix} (\emph{BDC matrix} for short) if
\begin{equation}
\mathbf{A}=\left[ \renewcommand{\arraystretch}{0.85}
\begin{array}{cccc}
\mathbf{C}_{0,0} & \mathbf{C}_{0,1} & \ldots  & \mathbf{C}_{0,t-1} \\
\mathbf{C}_{1,0} & \mathbf{C}_{1,1} & \ldots  & \mathbf{C}_{1,t-1} \\
\vdots  & \vdots  & \vdots  & \vdots  \\
\mathbf{C}_{t-1,0} & \mathbf{C}_{t-1,1} & \ldots  & \mathbf{C}_{t-1,t-1}
\end{array}
\right] ,\text{ }  \label{eq3.1 _block-circulant}
\end{equation}
where $\mathbf{C}_{i,j}$ is a \emph{circulant }$d\times d$
binary matrix for all $i,j$, and submatrices $\mathbf{C}_{i,j}$
and $\mathbf{C}_{l,m}$ with $ j-i\equiv m-l\pmod t$ have the
same weight. The matrix
\begin{equation}
\mathbf{W(A)}=\left[ \renewcommand{\arraystretch}{0.80}
\begin{array}{ccccccc}
w_{0} & w_{1} & w_{2} & w_{3} & \ldots  & w_{t-2} & w_{t-1} \\
w_{t-1} & w_{0} & w_{1} & w_{2} & \ldots  & w_{t-3} & w_{t-2} \\
w_{t-2} & w_{t-1} & w_{0} & w_{1} & \ldots  & w_{t-4} & w_{t-3} \\
\vdots  & \vdots  & \vdots  & \vdots  & \vdots  & \vdots  & \vdots  \\
w_{1} & w_{2} & w_{3} & w_{4} & \ldots  & w_{t-1} & w_{0}
\end{array}
\right]   \label{eq3.1_block-circulant-weights}
\end{equation}
is a \emph{circulant }$t\times t$ matrix whose entry in
position $i,j$ is the \emph{weight} of $\mathbf{C}_{i,j}$.
$\mathbf{W(A)}$ is called the \emph{weight matrix} of
$\mathbf{A.}$ The vector $\overline{\mathbf{W}}\mathbf{
(A)}=(w_{0},w_{1},\ldots ,w_{t-1}\mathbf{)}$ is called the
\emph{weight vector} of $\mathbf{A.}$
\end{definition}

\begin{remark}\label{zetamu} If in Definition 1 the matrices $C_{i,j}$ were assumed to be
right-circulant and not left-circulant, then they would have
been the sum of some right-circulant permutation matrices. A
right-circulant $d\times d$ permutation matrix is always
associated to a permutation of the set $\{1,2,\ldots,d\}$ in
the subgroup generated by the cycle $(1 \, 2 \, 3 \, \ldots\,
d)$. Then the notion of a BDC matrix is substantially
equivalent to that of a $\mathbb Z_d$-scheme, as defined in
\cite{AFLN-Semiplanes}.
\end{remark}

Let $\mathbf{A}$ be as in Definition \ref{def3.1_block double-circul}. In
addition, assume that $\mathbf{A}$ is the incidence matrix of a symmetric
configuration $v_{k}$ with $k=\sum_{i=0}^{t}w_{i}$. From $\mathbf{A}$ one
can obtain BDC incidence $v^{\prime }\times v^{\prime }$ matrices $\mathbf{A}
^{\prime }$ of symmetric configurations $v_{k^{\prime }}^{\prime }$ by the
following way.

(i) For each $h\in \{0,1,\ldots ,t-1\}$, in each row of every
submatrix $ \mathbf{C}_{i,j}$ with $j-i\equiv h\pmod t$ replace
$\delta _{h}\geq 0$ values of $1$ with zeros, in such a way
that the obtained submatrix is still circulant. As a result, a
BDC incidence matrix of a configuration $ v_{k^{\prime
}}^{\prime }$ with
\begin{equation}
v^{\prime }=v,\text{ }k^{\prime }=k-\sum_{h=0}^{t-1}\delta _{h},\text{ }
0\leq \delta _{h}\leq w_{h},\,w_{h}^{\prime }=w_{h}-\delta _{h},\overline{
\mathbf{W}}(\mathbf{A}^{\prime })=(w_{0}^{\prime },\ldots ,w_{t-1}^{\prime })
\label{eq3.1_rem(i)}
\end{equation}
is obtained.

(ii) Fix some non-negative integer $j\le t-1$. Let $m$ be such
that $ w_{m}\leq w_{h}$ for all $h\neq j$. Cyclically shift all
block rows of $ \mathbf{A}$ to the left by $j$ block positions.
A matrix $\mathbf{A}^{\ast } $ with
$\overline{\mathbf{W}}(\mathbf{A}^{\ast })=(w_{0}^{\ast
}=w_{j},\ldots ,w_{u}^{\ast }=w_{u+j\pmod t},\ldots
,w_{t-1}^{\ast }=w_{j-1}) $ is obtained. By applying (i),
construct a matrix $\mathbf{A}^{\ast \ast }$ with $w_{0}^{\ast
\ast }=w_{0}^{\ast }=w_{j},$ $w_{h}^{\ast \ast }=w_{m},$ $h\geq
1.$ Now remove from $\mathbf{A} ^{\ast \ast }$ $t-c$ block rows
and columns from the bottom and the right. In this way an
incidence $cd\times cd$ BDC matrix $\mathbf{A}^{\prime }$ of a
configuration $v_{k^{\prime }}^{\prime }$ is obtained with
\begin{equation}
v^{\prime }=cd,\,k^{\prime }=w_{j}+(c-1)w_{m},\,c=1,2,\ldots ,t,\,\overline{
\mathbf{W}}(\mathbf{A}^{\prime })=(w_{j},w_{m},\ldots ,w_{m}).
\label{eq3.1_rem(ii)}
\end{equation}

(iii) Let $t$ be even. Let $\mathbf{A}^{\ast }$ be as in (ii).
Let $w_{\text{O}},w_{\text{E}}$ be weights such that
$w_{\text{O}}\leq w_{h}^{\ast }$ for odd $h$ and
$w_{\text{E}}\leq w_{h}^{\ast }$ for even $h$. By applying (i),
construct a matrix $\mathbf{A}^{\ast \ast }$ with $w_{0}^{\ast
\ast }=w_{0}^{\ast }=w_{j},$ $w_{h}^{\ast \ast }=w_{\text{O}}$
for odd $h$, $ w_{h}^{\ast \ast }=w_{\text{E}}$ for even $h\geq
2$. Let $f=1,2,\ldots ,t/2$ . From $\mathbf{A}^{\ast \ast }$
remove $t-2f$ block rows and columns from the bottom and the
right. An incidence $2fd\times 2fd$ BDC matrix $\mathbf{A}
^{\prime }$ of a configuration $v_{k^{\prime }}^{\prime }$ is
obtained with
\begin{equation}
v^{\prime }=2fd,k^{\prime }=w_{j}+w_{\text{O}}+(f-1)(w_{\text{E}}+w_{\text{O
}}),\overline{\mathbf{W}}(\mathbf{A}^{\prime })=(w_{j},w_{\text{O}},
\underbrace{w_{\text{E}},w_{\text{O}},\ldots ,w_{\text{E}},w_{\text{O}}}_{f-1
\text{ pairs}}).  \label{eq3.1_rem(iii)}
\end{equation}

\begin{remark}\label{zetamu2}
Construction (i) in this section essentially follows from
Theorem \ref{Th2-Mend}; in \cite{AFLN-Semiplanes} it is
referred to as $1$-factor deletion. Construction (ii) is
essentialy a different formulation  of Proposition 4.3 in
\cite{AFLN-Semiplanes}, which is stated in terms of $\mathcal
S_\mu$-schemes. Apart from terminology, the only difference is
that here the case $m>1$ is considered. Other methods for
obtaining families of symmetric configurations from
$\mathbf{A}$  can be found in \cite[Sec.\thinspace
4]{DGMP-GraphCodes}.
\end{remark}

\section{Constructions and parameters of block \newline
double-circulant incidence matrices from
\cite{DFGMP-submitted}\label{sec_BDC}}

All the results of this section are taken from
\cite{DFGMP-submitted}, apart from Tables 3.1 and 3.2 which
essentially present more examples than the corresponding
\cite[Tab. 1]{DFGMP-submitted}. In Subsections \ref{sub_PG} and
\ref{sub_AG} we give some results based on a general method
connected with the action of the automorphism group $S$ of a
configuration. The method was originally proposed in
\cite{DGMP-ACCT2008,DGMP-GraphCodes} and then developed in
\cite{DFGMP-submitted}.

\subsection{BDC incidence matrices from projective planes\label{sub_PG}}

In this subsection, the projective plane $PG(2,q)$ is
considered as a cyclic symmetric configuration
$(q^{2}+q+1)_{q+1}$ \cite[ Sec.\thinspace 5]{DGMP-GraphCodes},
\cite[Sec.\thinspace 5.5]{Dimit}, \cite[ Th.\thinspace
19.15]{Shearer-Handb}, \cite{Singer}. The Singer group of
$PG(2,q)$ is used as the automorphism group $S.$

The following BDC\ matrices with $d\times d$ circulant submatrices and the
corresponding BDC\ configurations $v_{k}$ are given in \cite[Sec.\thinspace
4.1]{DFGMP-submitted}:
\begin{align}
\text{BDC }v_{k}& :\text{ }d=\frac{q^{2}+q+1}{3},\text{ }v=2d,\text{ }k=
\frac{2q+\sqrt{q}+2}{3},\text{ }q=p^{4m+2},\text{ }p\equiv 2\pmod 3; \\
\text{BDC }v_{k}& :\text{ }d=\frac{q^{2}+q+1}{3},\text{ }v=2d,\text{ }k=
\frac{2q-\sqrt{q}+2}{3},\text{ }q=p^{4m},\text{ }p\equiv 2\pmod 3.  \notag
\end{align}
\begin{align}
\text{BDC }v_{k}& :\text{ }d=\frac{q^{2}+q+1}{t},\text{ }v=cd,\text{ }k=
\frac{q+1\pm (1-t)\sqrt{q}}{t}+(c-1)\frac{q+1\pm \sqrt{q}}{t},
\label{eq3.3_Ex(ii)ProjPlane} \\
& c=1,2,\ldots ,t,\text{ }q=p^{2m},\text{ }t\text{ prime, }  \notag
\end{align}
where $p\pmod t$ is a generator of the multiplicative group of
${\mathbb{Z}} _{t}$.

The needed for \eqref{eq3.3_Ex(ii)ProjPlane} hypothesis that
$p\pmod t$ is a generator of the multiplicative group of
${\mathbb{Z}} _{t}$ holds, for example, in the following cases:
$q=3^{4},$ $ t=7;$ $q=2^{8},$ $t=13;$ $q=5^{4},$ $t=7;$
$q=2^{12},$ $t=19;$ $q=3^{8},$ $ t=7;$ $q=2^{16},$ $t=13;$
$q=17^{4},$ $t=7;$ $p\equiv 2\pmod t,$ $t=3.$

In Table 4.1, parameters of configurations $v_{n}^{\prime }$
with BDC incidence matrices are given. We use both (ii) and
(iii) of Section \ref{sec_bdc&families}. The starting weights
$w_{i}^{\ast }$ are obtained by computer by considering orbits
of subgroups of a Singer group of $PG(2,q)$. For $ q=81$ we
use~(\ref{eq3.3_Ex(ii)ProjPlane}). The values $k^{\prime
},v^{\prime }$ are calculated by
(\ref{eq3.1_rem(ii)}),(\ref{eq3.1_rem(iii)}). Only cases with
$v^{\prime }<G(k^{\prime })$ are included in the tables. Then
the smallest value $k^{\#}$ for which $v^{\prime }<G(k^{\#})$
is found. As a result, each row of the table provides
configurations $v_{n}^{\prime }$ with $v^{\prime }<G(n),$
$n=k^{\#},k^{\#}+1,\ldots ,k^{\prime },$ see (i) of Section
\ref{sec_bdc&families} and (\ref{eq3.1_rem(i)}). In column $
w_{i}^{\ast },$ an entry $s_{j}$ indicates that the
 weight $s$
should be repeated $j$ times.

\subsection{BDC incidence matrices from punctured affine planes\label{sub_AG}}

In this subsection, the cyclic punctured affine plane is
considered as a cyclic symmetric configuration $(q^{2}-1)_{q},$
see \cite{Bose},\cite[ Sec.\thinspace
5.6]{Dimit},\cite[Th.\thinspace 19.17]{Shearer-Handb}, as well
as \cite[Ex.\thinspace 5, Sec.\thinspace 6]{DGMP-GraphCodes}.
The affine Singer group of $AG(2,q)$ is used as the
automorphism group $S$.

The following BDC\ matrices with $d\times d$ circulant submatrices and the
corresponding BDC\ configurations $v_{k}$ are given in \cite[Sec.\thinspace
4.2]{DFGMP-submitted} on the base \cite[Th.\thinspace 4]{DFGMP-submitted}:
\begin{align}
\text{BDC }v_{k}& :\text{ }d=(\sqrt{q}-1)(q+1),\text{ }v=2fd,\text{ }
k=(2f-1) \sqrt{q},  \label{eq3.3_AfPlnPairsWeitghts} \\
\text{ }f& =1,2,\ldots ,\frac{\sqrt{q}+1}{2},\text{ }q\text{ odd square.}
\notag \\
\text{BDC }v_{k}& :\text{ }d=2(\sqrt{q}-1)(q+1),\text{ }v=cd,\text{ }
k=(2c-1) \sqrt{q},  \notag \\
\text{ }c& =1,2,\ldots ,\frac{\sqrt{q}+1}{2},\text{ }q\text{ odd square, }
\sqrt{q}\equiv 1\pmod 4.  \notag \\
\text{BDC }v_{k}& :\text{ }d=2(\sqrt{q}-1)(q+1),\text{ }v=2fd,\text{ }
k=(4f-1)\sqrt{q},  \notag \\
\text{ }f& =1,2,\ldots ,\frac{\sqrt{q}+1}{4},\text{ }q\text{ odd square, }
\sqrt{q}\equiv 3\pmod 4.  \notag \\
\text{BDC }v_{k}& :\text{ }d=4(\sqrt{q}-1)(q+1),\text{ }v=cd,\text{ }
k=(4c-1) \sqrt{q},  \notag \\
\text{ }c& =1,2,\ldots ,\frac{\sqrt{q}+1}{4},\text{ }q\text{ odd square, }
\sqrt{q}\equiv 3\pmod 4,\text{ }\frac{\sqrt{q}+1}{4}\text{ is odd.}  \notag
\\
\text{BDC }v_{k}& :\text{ }d=4(\sqrt{q}-1)(q+1),\text{ }v=2fd,\text{ }
k=(8f-1)\sqrt{q},  \notag \\
\text{ }f& =1,2,\ldots ,\frac{\sqrt{q}+1}{8},\text{ }q\text{ odd square, }
\sqrt{q}\equiv 3\pmod 4,\text{ }\frac{\sqrt{q}+1}{4}\text{ is even.}  \notag
\end{align}

In Table 4.2 parameters of configurations $v_{n}^{\prime }$
with BDC incidence matrices are given. We use both (ii) and
(iii) of Section \ref {sec_bdc&families}. The starting weights
$w_{i}^{\ast }$ are obtained by computer through the
constructions of the orbits of subgroups of the affine Singer
group. For notations $k^{\prime }$ and $k^{\#}$ see Table 4.1.

\section{Parameters of symmetric configurations $v_{k}$ admitting an
extension\label{sec_admitExten}}

The following infinite family of symmetric configuration
$v_{k}$ is obtained in \cite[Th. 6.2]{AFLN-Semiplanes},
\cite[Th. 1(i)]{AfDaZ-InfProc} with the help of
  Construction GM:
\begin{equation}
v_{k}:v=q^{2}-qs+\theta ,\text{ }k=q-s-\Delta ,\text{ }q>s\geq
0,\,q-s>\Delta \geq 0,\,\theta =0,1,\ldots ,q-s+1.  \label{eq2_tetaExten}
\end{equation}

\begin{theorem}
\emph{\cite[Corollary 1]{DFGMP-submitted}} Let $v=td$, $t\geq
k$, $d\geq k-1$ , and let $v_{k}$ be a symmetric configuration.
Assume that an incidence matrix $\mathbf{A}$ of $v_{k}$ is a
BDC matrix as in \emph{(\ref{eq3.1 _block-circulant})} with
weight vector ${\overline{\mathbf{W}}}(\mathbf{A})$.
If\thinspace\ $\overline{\mathbf{W}}(\mathbf{A})=(0,1,1,\ldots
,1)$ or\thinspace\
$\overline{\mathbf{W}}(\mathbf{A})=(1,1,\ldots ,1)$ then one
can obtain a family of symmetric configurations $v_{k}$ with
parameters \eqref{eq4_(011...1)} or \eqref{eq4_(11...1)},
respectively
\begin{eqnarray}
v_{k} &:&v=cd+\theta ,\text{ }k=c-1-\delta ,\text{ }c=2,3,\ldots ,t,\text{ }
\theta =0,1,\ldots ,c+1,\text{ }\delta \geq 0\text{.}  \label{eq4_(011...1)}
\\
v_{k} &:&v=cd+\theta ,\text{ }k=c-\delta ,\text{ }c=2,3,\ldots ,t,\text{ }
\theta =0,1,\ldots ,c+1,\text{ }\delta \geq 0\text{.}  \label{eq4_(11...1)}
\end{eqnarray}
\end{theorem}

Configurations with parameters \eqref{eq4_AfPlExten} were
first obtained in \cite[Th. 6.3]{AFLN-Semiplanes} from the punctured
affine plane by using Construction GM; see also \cite[Eqn.
(8)]{DGMP-Petersb2009} and \cite[Ex.~6(i)]{DFGMP-submitted}.
\begin{equation}
v_{k}:v=c(q-1)+\theta ,\text{ }k=c-1-\delta ,\text{ }c=2,3,\ldots ,q+1,\text{
}\theta =0,1,\ldots ,c+1,\text{ }\delta \geq 0\text{.}  \label{eq4_AfPlExten}
\end{equation}

The families of configurations $v_{k}$ (\ref{eq4_RuzExten}) and
(\ref{eq4_RuzExten_d=p-1}) are obtained in \cite[Sec.\thinspace
5]{DFGMP-submitted} by Construction GM, starting from some new
starting matrices proposed in \cite{DFGMP-submitted}.
\begin{equation}
v_{k}:v=cp+\theta ,\text{ }k=c-\delta ,\text{ }c=2,3,\ldots ,p-1,\text{ }
\theta =0,1,\ldots ,c+1,\text{ }\delta \geq 0,\text{ }p\text{ prime.}
\label{eq4_RuzExten}
\end{equation}
\begin{equation}
v_{k}:v=c(p-1)+\theta ,\text{ }k=c-1-\delta ,\text{ }c=2,3,\ldots ,p,\text{ }
\theta =0,1,\ldots ,c+1,\,\delta \geq 0,\,p\text{ prime.}
\label{eq4_RuzExten_d=p-1}
\end{equation}
The family of configurations $v_{k}$ (\ref{eq4_BaerExten}) is
obtained in \cite{DFGMP-submitted}.
\begin{eqnarray}
v_{k} &:&v=c(q+\sqrt{q}+1)+\theta ,\text{ }k=c-\delta ,\text{ }c=2,3,\ldots
,q-\sqrt{q}+1,  \notag \\
&&\theta =0,1,\ldots ,c+1,\text{ }\delta \geq 0\text{, }q\text{ square}.
\label{eq4_BaerExten}
\end{eqnarray}

\section{The spectrum of parameters of cyclic symmetric configurations \label{sec_ParamCyclicConfig}}

Current data on the existence of cyclic configurations $v_{k}$,
$k\leq 51$, are given in Table~6.1. New parameters obtained in
\cite{DFGMP-submitted} are written in bold font.

In Table 6.1, the values of $v$ for which cyclic symmetric
configurations $ v_{k}$ exist (resp. do not exist) are written
in normal (resp. in italic) font. Moreover, $\overline{v}$
means that no configuration $v_{k}$ exists, whereas
$\overline{v^{c}}$ indicates that no cyclic configuration
$v_{k}$ exists. Data from
\cite{Funk2008,Gropp-nk,Gropp-nonsim,Gropp-Handb,KaskiOst,Lipman,ShearerWebModulGR}
are used in the 4-th column of the table. We take into account
that an entry of the form \textquotedblleft
$t+$\textquotedblright\ in the row \textquotedblleft
$n$\textquotedblright\ of \cite[Tab.\thinspace
1]{ShearerWebModulGR} means the existence of cyclic symmetric
configurations $ v_{n}$ with $v\geq t.$ An absence of a value
$"v"$ in the row \textquotedblleft $n$\textquotedblright\ of
\cite[Tab.\thinspace 1]{ShearerWebModulGR} means the
non-existence of a cyclic symmetric configuration $v_{n}$.
Also, we use the following \emph{non-existence} results:
$\overline{\emph{32} _{6}}$ \cite[Th.\thinspace 4.8]{Gropp-nk};
$\overline{\emph{33}_{6}}$ \cite {KaskiOst};
$\overline{\emph{34}_{6}^{c}},$
$\overline{\emph{59}_{8}^{c}}$-$ \overline{\emph{62}_{8}^{c}}$
\cite{Lipman}; $\overline{\emph{75}_{9}^{c}}$-$
\overline{\emph{79}_{9}^{c}},\overline{\emph{81}_{9}^{c}}$-
$\overline{\emph{84}^{c}_{9}}$ \cite{Funk2008};
$\overline{\emph{93}^{c}_{10}}-\overline{\emph{106}^{c}_{10}}$,
$\overline{\emph{121}^{c}_{11}}-\overline{\emph{132}^{c}_{11}}$,
$\overline{\emph{134}^{c}_{11}},$
$\overline{\emph{135}^{c}_{12}}-\overline{\emph{155}^{c}_{12}}$,
$\overline{\emph{157}^{c}_{12}}$,
$\overline{\emph{160}^{c}_{12}}$,
$\overline{\emph{169}^{c}_{13}}-\overline{\emph{182}^{c}_{13}}$,
$\overline{\emph{184}^{c}_{13}}-\overline{\emph{192}^{c}_{13}}$,
$\overline{\emph{185}^{c}_{14}}-\overline{\emph{224}^{c}_{14}}$,
$\overline{\emph{256}^{c}_{15}}-\overline{\emph{260}^{c}_{15}}$,
$\overline{\emph{263}^{c}_{15}}$ \cite[Tab.\thinspace
 1]{ShearerWebModulGR}; see also Theorem
\ref{th5_deficiency1}.

The values of $k$ for which the spectrum of parameters of
cyclic symmetric configurations $v_{k}$ is completely known are
indicated by a dot "$ \centerdot $"; the corresponding values
of $E_{c}(k)$ are sharp and they are noted by the dot
"$\centerdot $" too.

An entry $\mathbf{v}$-$\mathbf{w}$ indicates an interval of sizes
from $\mathbf{v}$ to $\mathbf{w}$ without gaps. If an already
known value lies within an interval $\mathbf{v}$-$\mathbf{w}$
obtained in work \cite{DFGMP-submitted}, then it is written
immediately before the interval.

Entries $v_{a},v_{b},$ and $v_{c}$ (here and in all tables)
mean, respectively, that relations (\ref{eq2_cyclicPG(2,q)}),
(\ref{eq2_cyclicAG(2,q)}), and (\ref{eq2_cyclicRuzsa}) are
applied.

The value $v_{\delta }(k)$ is defined in Introduction. In the
second column, the exact values of $v_{\delta }(k)$ are marked
by the dot \textquotedblleft $\centerdot $\textquotedblright .
For $k\leq 16$, the exact values of $ v_{\delta }(k)$ are taken
from \cite[Tab.\thinspace IV]{GrahSloan}, \cite[ Tab.\thinspace
2]{OstergModulSidon}, \cite[Tab.\thinspace
1a]{ShearerWebModulGR}, \cite{Swanson}. Also, if $k-1$ is a
prime power then $ v_{\delta }(k)=P(k)=k^{2}-k+1$. The
remaining entries in the second column are lower bounds of
$v_{\delta }(k)$. By the Bruck-Ryser Theorem, planes $
PG(2,k-1)$ with $k-1=6,14,21,22,30,33,38,42,46,54,57,62$ do not
exist. It is well-known that $P(k)\leq v_{\delta }(k)$, and
that a \emph{cyclic} symmetric configuration $(k^{2}-k+1)_{k}$
exists if and only if a \emph{cyclic} projective plane of order
$k-1$ exists. By \cite{BaumertGordon}, no cyclic projective
planes exist with non-prime power orders $\leq 2\cdot 10^{9}$.
Therefore cyclic projective planes $PG(2,k-1)$ with $
k-1=18,20,24,26,28,34,35,36,39,40,44,45,48,50$ do not exist.
Also we use Theorem~\ref{th5_deficiency1}
taken from~\cite{Gropp-nonsim}. The mentioned
 non-existence cases of
cyclic configurations are marked in Table~6.1 by subscripts
$br$ (Bruck-Ryser Theorem), $s$ (\cite{BaumertGordon}), and
$t$ (Theorem~\ref{th5_deficiency1}).

For $k\le22$ the \emph{filling of the interval} $P(k)-G(k)$ is expressed as
a percentage in the last column of Table 6.1.

\begin{theorem}
\label{th5_deficiency1}\emph{\cite[Th.\thinspace 2.4]{Gropp-nonsim}} There
is no symmetric configuration $(k^{2}-k+2)_{k}$ if $\,5\leq k\leq 10$ or if
neither $k$ or $k-2$ is a square.
\end{theorem}

In order to widen the ranges of parameter pairs
 $(v,k)$ for which a cyclic
symmetric configuration $v_{k}$ exists, we consider a number of
procedures that allow to define a new modular Golomb ruler from
a known one. Some methods have already been introduced in the
paper, see Theorem \ref{Th2-Mend}.

Here we first recall a result from \cite{Shearer-Handb}, which describes a
method to construct different rulers with the same parameters.

\begin{theorem}
\label{Th2_newMGRfromold} \emph{\cite{Shearer-Handb}} If $
(a_{1},a_{2},\ldots ,a_{k})$ is a $(v,k)$ modular Golomb ruler
and $m$ and $ b $ are integers with $\gcd (m,v)=1$ then
$(ma_{1}+b\pmod v,ma_{2}+b\pmod v ,\ldots ,ma_{k}+b\pmod v)$ is
also a $(v,k)$ modular Golomb ruler.
\end{theorem}

It should be noted that a $(v,k)$ modular Golomb ruler can be a
$(v+\Delta ,k)$ modular Golomb ruler for some integer $\Delta $
\cite{Gropp-nk}. This property does not depend on parameters
$v$ and $k$ only. This is why Theorem \ref{Th2_newMGRfromold}
can be useful for our purposes.

\begin{example}
\label{ex5_extenGR} We consider the $(31,6)$ modular Golomb ruler
\begin{equation*}
(a_{1},\ldots ,a_{6})=(0,1,4,10,12,17)
\end{equation*}
obtained from $PG(2,5),$ see \cite{ShearerWebShortest}. We can
apply Theorem \ref{Th2_newMGRfromold} for $m=19$, $b=0$. The
$(31,6)$ modular Golomb ruler $(ma_{1}$ $(\text{mod }$
$31),\ldots ,ma_{6}$ $(\text{mod }31))$ is
\begin{equation*}
(a_{1}^{\prime },\ldots ,a_{6}^{\prime })=(0,4,11,13,14,19).
\end{equation*}
Now we take $\Delta =4$ and calculate the set of differences $
\{a_{i}^{\prime }-a_{j}^{\prime }\,(\text{mod }35)|\,1\leq
i,j\leq 6;i\neq j\}$, that is $
\{1,2,3,4,5,6,7,8,9,10,11,13,14,15,16,19,20,21,22,24,25,26,27,28,29,30,
\linebreak 31,32,33,34\}.$ As the all differences are distinct
and nonzero, the starting $(31,6)$ modular Golomb ruler is also
a $(35,6)$ modular Golomb ruler.
\end{example}

For $k\leq 81$, we performed a computer search starting from
the $(v,k)$ modular Golomb rulers corresponding to
(\ref{eq2_cyclicPG(2,q)})--(\ref{eq2_cyclicRuzsa}). For
projective and affine planes, we got a concrete description of
the ruler from \cite{ShearerWebShortest}.

For Ruzsa's construction, we use the following known relations. Let $p$ be a
prime. Let $g$ be a primitive element of $F_{p}$. The following Ruzsa's
sequence \cite{Ruzsa},\cite[Sec.\thinspace 5.4]{Dimit},\cite[Th.\thinspace
19.19]{Shearer-Handb} forms a $(p^{2}-p,p-1)$ modular Golomb ruler:
\begin{equation}
e_{u}=pu+(p-1)g^{u}\pmod{p^2-p},\text{ }u=1,2,\ldots ,p-1,\text{ }v=p^{2}-p.
\label{eq3.2_Ruzsa}
\end{equation}

For every starting $(v,k)$ modular ruler we first considered
all possible $m$ with
\begin{center}
$\gcd (m,v)=1$,
\end{center}
and applied Theorem \ref{Th2_newMGRfromold} for all $b<v$ to
get new rulers with the same parameters $v$ and~$k$. Then, we
checked whether this ruler was also a $(v+\Delta ,k)$ for some
$\Delta $.

In Table 6.2, for $52\leq k\leq 83,$ the upper bounds on the
cyclic existence bound $E_{c}(k)$ obtained in
\cite{DFGMP-submitted} are listed in bold font. An
 entry, say
$A(u)$, in the column $E_{c}(k)$ on the row \textquotedblleft
$u$\textquotedblright, means that in \cite{DFGMP-submitted}
\emph{all} cyclic symmetric configurations $v_{u}$ in the
region $A(u),A(u)+1,\ldots ,G(u)-1$ are obtained. These
configurations are new. We obtained also many other new cyclic
symmetric configurations $ v_{k}$ for $52\leq k\leq 83$.
However, we do not give here their sizes $v$ here in order to
save space.

The upper bounds on $E_{c}(k)$ in Tables 6.1 and 6.2, obtained
in\cite{DFGMP-submitted}, are written in bold font.

Some new cyclic symmetric configurations important for Table
7.1 of Section 7 are given in Table 6.3 where we write the
first rows of their incidence matrices; these rows may be the
same for distinct $v$.

\section{The spectrum of parameters of symmetric (not necessarily cyclic)
configurations \label{sec6_spectrum}}

The known results regarding to parameters of symmetric
configurations can be found in
\cite{AFLN-graphs,AFLN-ConfigGraphs,AFLN-Semiplanes,AFLN-CyclSchem,AfDaZ,AfDaZ-InfProc,%
ArParBalbNetw2011,ArParBalbHegDM2010,Baker,BalbuenaSIAM2008,Boben-v3,Bose,%
CarsDinStef-Reduc-n3,DFGMPConfigArxiv2012,DFGMP-submitted,DGMP-ACCT2008,%
DGMP-Petersb2009,DGMP-GraphCodes,Funk1993,Funk2008,FunkLabNap,GH,GrahSloan,Gropp-nk,Gropp-Chemic,Gropp-nonsim,Gropp-ConfGraph,%
Gropp-ConfGeomCombin,Gropp-Handb,Grunbaum,OstergModulSidon,KaskiOst,Krcadinac,%
Lipman,Martinetti,MePaWolk,Ruzsa,Shearer-Handb,ShearerWebModulGR};
see also the references therein.

The known families of configurations are described in Section
\ref{sec_known}, see also \eqref{eq2_tetaExten},\eqref{eq4_AfPlExten}. New
families obtained in the work \cite{DFGMP-submitted} are given
in Sections \ref{sec_BDC} and \ref{sec_admitExten}. In Table
7.1, for $k\leq 51,$ $P(k)\leq v<G(k),$ values of $v$ for which
a symmetric configuration $v_{k}$ from one of the families of
Sections \ref{sec_known}--\ref{sec_admitExten} exists are
given. The new parameters obtained in the paper
\cite{DFGMP-submitted} are written in bold font.

An entry of type $v_{\text{subscript}}$ indicates that one of
the following is used: (2.$i$), (4.$ j$), (5.$k$), Table
4.1,Table 6.1, the Bruck-Ryser Theorem, Theorem
\ref{th5_deficiency1}. More precisely $v_{a}$ indicates that
$v$ is obtained from $(\ref{eq2_cyclicPG(2,q)})$, and similarly
$v_{b}\rightarrow ( \ref{eq2_cyclicAG(2,q)}),$
$v_{c}\rightarrow (\ref{eq2_cyclicRuzsa}),$ $ v_{f}\rightarrow
(\ref{eq2_q-1-cancel}),$ $v_{g}\rightarrow (\ref{eq2_Baer}), $
$v_{h}\rightarrow (\ref{eq2_Baer2}),$ $v_{j}\rightarrow (\ref
{eq2_FuLabNabDecomp}),$ $v_{k}\rightarrow
(\ref{eq2_affine_q-1}),$ $v_{\lambda}\rightarrow
\eqref{eq2_Balbuena1}-\eqref{eq2_AaParBalb}$, $
v_{m}\rightarrow (\ref{eq2_tetaExten})$, $v_{P}\rightarrow
(\ref{eq3.3_AfPlnPairsWeitghts})$, $v_{r}\rightarrow
(\ref{eq4_AfPlExten})$, $ v_{S}\rightarrow
(\ref{eq4_RuzExten})$, $v_{T}\rightarrow
(\ref{eq4_RuzExten_d=p-1})$, $v_{W}\rightarrow $ Table~4.1,
$v_{y}\rightarrow $ Table~6.1 with $k\leq 15$,
$v_{Z}\rightarrow $ Table~6.1 with $k>15$, $ v_{br}\rightarrow
$ the Bruck-Ryser Theorem, $v_{t}\rightarrow $Theorem
\ref{th5_deficiency1}. Here capital letters in subscripts
remark new results and constructions of \cite{DFGMP-submitted}
while lower case letters indicate the known ones.

An entry $v_{\text{subscript}_{1}\cdot
\text{subscript}_{2}\cdot \ldots }$ with more than one
subscript means that the same value can be obtained from
different constructions. An entry of type
$v_{\text{subscript}_{1}\cdot \text{subscript}_{2}\cdot \ldots
}-v_{\text{subscript}_{1}\cdot \text{ subscript}_{2}\cdot
\ldots }^{\prime }$ indicates that a whole interval of values
from $v$ to $v^{\prime }$ can be obtained from the
constructions corresponding to the subscripts. We use the
following known results on the existence of sporadic symmetric
configurations: $45_{7}$ \cite{Baker}; $ 82_{9}$
\cite[Tab.\thinspace 1]{FunkLabNap}; $135_{12},$ see
\cite{Gropp-nonsim} with reference to Mathon's talk at the
British Combinatorial Conference 1987; $34_{6}$
\cite{Krcadinac}, see also \cite{AFLN-CyclSchem}. The
non-existence of configuration $112_{11}$ is proven in
\cite{KaskOst112-11}. The non-existence of the plane $PG(2,10)$
implies the non-existence of configuration $111_{11}$. The
values of $k$ for which the spectrum of parameters of symmetric
configurations $v_{k}$ is completely known are indicated by a
dot "$\centerdot $"; the corresponding values of $E(k)$ are
exact and they are indicated by a dot as well.

To save space, in Table 7.1 for
 given $v,k$, we do not write all
the constructions providing a configuration $v_{k}$, but often we describe some of them
as a matter of illustration.

For the convenience of the reader we give also Table 7.2 where
constructions are not indicated and for $k\leq 64,$ $P(k)\leq
v<G(k),$ values of $v$ for which a symmetric configuration
$v_{k}$ exists are written.

The {filling of the interval} $P(k)-G(k)$ is expressed as a
percentage in the last column of Tables 7.1 and 7.2. It is
interesting to note that such a percentage is quite high and
that most gaps occur for $v$ close to $k^{2}-k+1 $.

We note that a number of parameters obtained in the work
\cite{DFGMP-submitted} are new, see bold font in Table 7.1.
Note that parameters of some new families are too big to be
included in Tables 7.1 and 7.2. Recall also (see Introduction)
that from the stand point of applications, including Coding
Theory, it is useful to have different matrices
${\mathbf{M}}(v,k)$ for the same $v$ and $k$.


\section{Tables of parameters \label{sec_tables}}

\subsection{Tables for Section \protect\ref{sec_BDC}}

\begin{center}
\textbf{Table 4.1. }Parameters of configurations $v_{n}^{\prime }$ with BDC
incidence matrices$\mathbf{,}$ $v^{\prime }<G(n),$ $n=k^{\#},k^{\#}+1,\ldots
,k^{\prime },$ by (ii) of Section \ref{sec_bdc&families} from the cyclic
projective plane $PG(2,q)\smallskip $

$
\begin{array}{r|c|r|c|r|r|r|r|r|r}
\hline
q & t & d & w_{i}^{\ast } & c & k^{\prime } & v^{\prime } & G(k^{\prime }) &
k^{\#} & G(k^{\#}) \\ \hline
25 & 3 & 217 & 12,7,7 & 2 & 19 & 434 & 493 & 19 & 493 \\ \hline
25 & 7 & 93 & 8,3,3,3,3,3,3 & 4 & 17 & 372 & 399 & 17 & 399 \\ \hline
25 & 7 & 93 & 8,3,3,3,3,3,3 & 5 & 20 & 465 & 567 & 19 & 493 \\ \hline
25 & 7 & 93 & 8,3,3,3,3,3,3 & 6 & 23 & 558 & 745 & 20 & 567 \\ \hline
32 & 7 & 151 & 0,5,5,6,5,6,6 & 6 & 25 & 906 & 961 & 25 & 961 \\ \hline
37 & 3 & 469 & 16,9,13 & 2 & 25 & 938 & 961 & 25 & 961 \\ \hline
43 & 3 & 631 & 19,13,12 & 2 & 31 & 1262 & 1495 & 30 & 1361 \\ \hline
49 & 3 & 817 & 21,16,13 & 2 & 34 & 1634 & 1877 & 33 & 1719 \\ \hline
61 & 3 & 1261 & 25,21,16 & 2 & 41 & 2522 & 2611 & 40 & 2565 \\ \hline
64 & 3 & 1387 & 27,19,19 & 2 & 46 & 2774 & 3407 & 42 & 2795 \\ \hline
64 & 19 & 219 & 11,3_{18} & 7 & 32 & 1533 & 1569 & 32 & 1569 \\ \hline
64 & 19 & 219 & 11,3_{18} & 8 & 35 & 1752 & 1975 & 34 & 1877 \\ \hline
64 & 19 & 219 & 11,3_{18} & 9 & 38 & 1971 & 2293 & 35 & 1975 \\ \hline
64 & 19 & 219 & 11,3_{18} & 10 & 41 & 2190 & 2611 & 37 & 2199 \\ \hline
64 & 19 & 219 & 11,3_{18} & 17 & 62 & 3723 & 6431 & 48 & 3775 \\ \hline
64 & 19 & 219 & 11,3_{18} & 18 & 65 & 3942 & 7187 & 50 & 4189 \\ \hline
67 & 3 & 1519 & 28,19,21 & 2 & 47 & 3038 & 3609 & 44 & 3193 \\ \hline
73 & 3 & 1801 & 28,27,19 & 2 & 47 & 3602 & 3609 & 47 & 3609 \\ \hline
79 & 3 & 2107 & 31,21,28 & 2 & 52 & 4214 & 4541 & 51 & 4381 \\ \hline
81 & 7 & 949 & 4,13,13,13,13,13,13 & 6 & 69 & 5694 & 8291 & 58 & 5703 \\
\hline
81 & 7 & 949 & 4,13,13,13,13,13,13 & 5 & 56 & 4745 & 5451 & 54 & 4747 \\
\hline
97 & 3 & 3169 & 39,28,31 & 2 & 67 & 6338 & 7639 & 62 & 6431 \\ \hline
103 & 3 & 3571 & 39,28,37 & 2 & 67 & 7142 & 7639 & 65 & 7187 \\ \hline
107 & 7 & 1651 & 24,15,15,13,15,13,13 & 6 & 89 & 9906 & 13557 & 75 & 9965 \\
\hline
107 & 7 & 1651 & 24,15,15,13,15,13,13 & 5 & 76 & 8255 & 10179 & 69 & 8291 \\
\hline
109 & 3 & 3997 & 43,36,31 & 2 & 74 & 7994 & 9507 & 69 & 8291 \\ \hline
109 & 7 & 1713 & 8,15,15,19,15,19,19 & 6 & 83 & 10278 & 12041 & 77 & 10409
\\ \hline
121 & 37 & 399 & 14,3_{36} & 25 & 89 & 9975 & 13557 & 76 & 10179 \\ \hline
127 & 3 & 5419 & 49,43,36 & 2 & 85 & 10838 & 12821 & 80 & 11127 \\ \hline
128 & 7 & 2359 & 24,21,21,14,21,14,14 & 6 & 94 & 14154 & 15769 & 91 & 15085
\\ \hline
137 & 7 & 2701 & 24,15,15,23,15,23,23 & 6 & 99 & 16206 & 17081 & 96 & 16243
\\ \hline
\end{array}
$\newpage

\textbf{Table 4.1 }(continue).\textbf{\ }Parameters of
configurations $ v_{n}^{\prime }$ with BDC incidence
matrices$\mathbf{,}$ $v^{\prime }<G(n),$
$n=k^{\#},k^{\#}+1,\ldots ,k^{\prime },$ by (ii) of Section
\ref{sec_bdc&families} from the cyclic projective plane
$PG(2,q)\smallskip $

$
\begin{array}{r|c|r|c|r|r|r|r|r|r}
\hline
q & t & d & w_{i}^{\ast } & c & k^{\prime } & v^{\prime } & G(k^{\prime }) &
k^{\#} & G(k^{\#}) \\ \hline
139 & 3 & 6487 & 52,39,49 & 2 & 91 & 12974 & 15085 & 86 & 13075 \\ \hline
149 & 7 & 3193 & 12,25,25,21,25,21,21 & 6 & 117 & 19158 & 25035 & 104 & 19163
\\ \hline
149 & 7 & 3193 & 12,25,25,21,25,21,21 & 5 & 96 & 15965 & 16243 & 96 & 16243
\\ \hline
121 & 3 & 4921 & 48,37,37 & 2 & 85 & 9842 & 12821 & 75 & 9965 \\ \hline
121 & 7 & 2109 & 21,20,13,13,21,13,21 & 6 & 86 & 12654 & 13075 & 85 & 12821
\\ \hline
121 & 37 & 399 & 14,3_{36} & 20 & 71 & 7980 & 8661 & 69 & 8291 \\ \hline
151 & 3 & 7651 & 57,43,52 & 2 & 100 & 15302 & 17663 & 93 & 15453 \\ \hline
151 & 7 & 3279 & 32,19,19,21,19,21,21 & 6 & 127 & 19674 & 28921 & 105 & 19769
\\ \hline
151 & 7 & 3279 & 32,19,19,21,19,21,21 & 5 & 108 & 16395 & 20831 & 97 & 16715
\\ \hline
151 & 7 & 3279 & 32,19,19,21,19,21,21 & 4 & 89 & 13116 & 13557 & 87 & 13417
\\ \hline
157 & 3 & 8269 & 61,48,49 & 2 & 109 & 16538 & 21167 & 97 & 16715 \\ \hline
163 & 3 & 8911 & 63,49,52 & 2 & 112 & 17822 & 27043 & 102 & 18437 \\ \hline
163 & 7 & 3819 & 32,25,25,19,25,19,19 & 6 & 127 & 22914 & 28921 & 114 & 23529
\\ \hline
163 & 7 & 3819 & 32,25,25,19,25,19,19 & 5 & 108 & 19095 & 20831 & 104 & 19163
\\ \hline
169 & 3 & 9577 & 64,57,49 & 2 & 113 & 19154 & 22847 & 104 & 19163 \\ \hline
179 & 7 & 4603 & 24,21,21,31,21,31,31 & 6 & 129 & 27618 & 30151 & 124 & 27897
\\ \hline
181 & 3 & 10981 & 67,63,52 & 2 & 119 & 21962 & 25823 & 111 & 22217 \\ \hline
\end{array}
$\newpage

\textbf{Table 4.2. }Parameters of configurations $v_{n}^{\prime }$ with BDC
incidence matrices$\mathbf{,}$ $v^{\prime }<G(n),$ $n=k^{\#},k^{\#}+1,\ldots
,k^{\prime },$ by (ii) and (iii) of Section \ref{sec_bdc&families} from the
cyclic punctured affine plane $AG(2,q)\smallskip $

$\renewcommand{\arraystretch}{1.0}
\begin{array}{r|c|r|c|c|r|r|r|r|r|r}
\hline
q & t & d & w_{i}^{\ast } & c & f & k^{\prime } & v^{\prime } & G(k^{\prime
}) & k^{\#} & G(k^{\#}) \\ \hline
16 & 3 & 85 & 8,4,4 & 2 &  & 12 & 170 & 171 & 12 & 171 \\ \hline
31 & 3 & 320 & 14,9,8 & 2 &  & 22 & 640 & 713 & 21 & 667 \\ \hline
37 & 3 & 456 & 16,9,12 & 2 &  & 25 & 912 & 961 & 25 & 961 \\ \hline
49 & 4 & 600 & 16,12,9,12 & 3 &  & 34 & 1800 & 1877 & 34 & 1877 \\ \hline
49 & 6 & 400 & 4,9,12,8,8,8 & 5 &  & 36 & 2000 & 2011 & 36 & 2011 \\ \hline
53 & 4 & 702 & 17,12,10,14 & 3 &  & 37 & 2106 & 2199 & 37 & 2199 \\ \hline
61 & 3 & 1240 & 25,16,20 & 2 &  & 41 & 2480 & 2611 & 39 & 2505 \\ \hline
61 & 4 & 930 & 18,18,12,13 & 3 &  & 42 & 2790 & 2795 & 42 & 2795 \\ \hline
61 & 5 & 744 & 18,9,10,12,12 & 4 &  & 45 & 2976 & 3375 & 43 & 3015 \\ \hline
61 & 6 & 620 & 15,8,8,10,8,12 & 5 &  & 47 & 3100 & 3609 & 44 & 3193 \\ \hline
64 & 9 & 455 & 0,8,8,8,8,8,8,8,8 & 8 &  & 56 & 3640 & 5451 & 48 & 3775 \\
\hline
64 & 9 & 455 & 0,8,8,8,8,8,8,8,8 & 7 &  & 48 & 3185 & 3775 & 44 & 3193 \\
\hline
67 & 3 & 1496 & 26,24,17 & 2 &  & 43 & 2992 & 3015 & 43 & 3015 \\ \hline
71 & 5 & 1008 & 8,14,15,16,18 & 4 &  & 50 & 4032 & 4189 & 50 & 4189 \\ \hline
73 & 3 & 1776 & 30,21,22 & 2 &  & 51 & 3552 & 4381 & 47 & 3609 \\ \hline
79 & 3 & 2080 & 32,25,22 & 2 &  & 54 & 4160 & 4747 & 50 & 4189 \\ \hline
79 & 6 & 1040 & 8,14,13,14,18,12 & 5 &  & 56 & 5200 & 5451 & 56 & 5451 \\
\hline
79 & 6 & 1040 & 18,12,8,14,13,14 &  & 2 & 50 & 4160 & 4189 & 50 & 4189 \\
\hline
81 & 4 & 1640 & 25,20,16,20 & 3 &  & 57 & 4920 & 5547 & 55 & 5197 \\ \hline
81 & 4 & 1640 & 25,20,16,20 &  & 1 & 45 & 3280 & 3375 & 45 & 3375 \\ \hline
81 & 8 & 820 & 16,8,8,8,9,12,8,12 & 7 &  & 64 & 5740 & 7055 & 59 & 5823 \\
\hline
81 & 8 & 820 & 16,8,8,8,9,12,8,12 & 6 &  & 56 & 4920 & 5451 & 55 & 5187 \\
\hline
83 & 8 & 861 & 6,10,10,10,15,11,10,11 & 7 &  & 66 & 6027 & 7515 & 60 & 6039
\\ \hline
83 & 8 & 861 & 6,10,10,10,15,11,10,11 & 6 &  & 56 & 5166 & 5451 & 55 & 5187
\\ \hline
89 & 4 & 1980 & 26,25,18,20 & 3 &  & 62 & 5940 & 6431 & 60 & 6039 \\ \hline
89 & 8 & 990 & 17,8,10,12,8,10,10,14 & 7 &  & 65 & 6930 & 7187 & 64 & 7055
\\ \hline
97 & 3 & 3136 & 37,34,26 & 2 &  & 63 & 6272 & 6783 & 62 & 6431 \\ \hline
97 & 4 & 2352 & 29,26,20,22 & 3 &  & 69 & 7056 & 8291 & 65 & 7187 \\ \hline
97 & 6 & 1586 & 21,20,14,16,14,12 & 5 &  & 69 & 7930 & 8291 & 69 & 8291 \\
\hline
101 & 4 & 2550 & 30,25,20,26 & 3 &  & 70 & 7650 & 8435 & 68 & 7913 \\ \hline
101 & 4 & 2550 & 30,25,20,26 &  & 1 & 55 & 5100 & 5197 & 55 & 5197 \\ \hline
103 & 3 & 3536 & 41,32,30 & 2 &  & 71 & 7072 & 8661 & 65 & 7187 \\ \hline
103 & 6 & 1768 & 24,18,14,17,14,16 & 5 &  & 80 & 8840 & 11127 & 72 & 8947 \\
\hline
103 & 6 & 1768 & 24,18,14,17,14,16 & 4 &  & 66 & 7072 & 7515 & 65 & 7187 \\
\hline
103 & 6 & 1768 & 24,18,14,17,14,16 &  & 2 & 70 & 7072 & 8435 & 65 & 7187 \\
\hline
107 & 8 & 1431 & 17,15,10,13,10,12,17,13 & 7 &  & 77 & 10017 & 10409 & 76 &
10179 \\ \hline
107 & 8 & 1431 & 17,15,10,13,10,12,17,13 &  & 3 & 73 & 8586 & 9027 & 71 &
8661 \\ \hline
\end{array}
$\newpage

\textbf{Table 4.2} (continue)\textbf{. }Parameters of
configurations $ v_{n}^{\prime }$ with BDC incidence
matrices$\mathbf{,}$ $v^{\prime }<G(n),$
$n=k^{\#},k^{\#}+1,\ldots ,k^{\prime },$ by (ii) and (iii) of
Section \ref{sec_bdc&families} from the cyclic punctured affine
plane $ AG(2,q)\smallskip $

$
\begin{array}{r|c|r|c|c|r|r|r|r|r|r}
\hline
q & t & d & w_{i}^{\ast } & c & f & k^{\prime } & v^{\prime } & G(k^{\prime
}) & k^{\#} & G(k^{\#}) \\ \hline
109 & 3 & 3960 & 42,30,37 & 2 &  & 72 & 7920 & 8947 & 69 & 8291 \\ \hline
109 & 4 & 2970 & 32,26,22,29 & 3 &  & 76 & 8910 & 10179 & 72 & 8947 \\ \hline
109 & 6 & 1980 & 12,19,24,18,18,18 & 5 &  & 84 & 9900 & 12319 & 75 & 9965 \\
\hline
109 & 9 & 1320 & 8,14,20,10,13,10,12,10,12 & 8 &  & 78 & 10560 & 10599 & 78
& 10599 \\ \hline
113 & 4 & 3192 & 32,32,24,25 & 3 &  & 80 & 9576 & 11127 & 75 & 9965 \\ \hline
113 & 7 & 1824 & 10,20,14,17,22,16,14 & 6 &  & 80 & 10944 & 11127 & 80 &
11127 \\ \hline
121 & 4 & 3660 & 36,30,25,30 & 3 &  & 86 & 10980 & 13075 & 80 & 11127 \\
\hline
121 & 4 & 3660 & 36,30,25,30 &  & 1 & 66 & 7320 & 7515 & 66 & 7515 \\ \hline
121 & 5 & 2928 & 16,24,28,25,28 & 4 &  & 88 & 11712 & 13491 & 83 & 12041 \\
\hline
125 & 4 & 3906 & 37,32,26,30 & 3 &  & 89 & 11718 & 13557 & 83 & 12041 \\
\hline
125 & 8 & 1953 & 24,16,14,15,13,16,12,15 & 7 &  & 96 & 13671 & 16243 & 90 &
13935 \\ \hline
125 & 8 & 1953 & 24,16,14,15,13,16,12,15 &  & 3 & 93 & 11718 & 15453 & 83 &
12041 \\ \hline
127 & 3 & 5376 & 49,36,42 & 2 &  & 85 & 10752 & 12821 & 79 & 10817 \\ \hline
127 & 6 & 2688 & 28,18,18,21,18,24 & 5 &  & 100 & 13440 & 17663 & 88 & 13491
\\ \hline
127 & 6 & 2688 & 28,18,18,21,18,24 & 4 &  & 82 & 10752 & 11629 & 79 & 10817
\\ \hline
131 & 5 & 3432 & 19,24,32,26,30 & 4 &  & 91 & 13728 & 15085 & 90 & 13935 \\
\hline
137 & 4 & 4692 & 40,32,29,36 & 3 &  & 98 & 14076 & 16925 & 91 & 15085 \\
\hline
139 & 3 & 6440 & 54,41,44 & 2 &  & 95 & 12880 & 15935 & 86 & 13075 \\ \hline
149 & 4 & 5550 & 42,41,32,34 & 3 &  & 106 & 16650 & 20271 & 97 & 16715 \\
\hline
151 & 3 & 7600 & 58,44,49 & 2 &  & 102 & 15200 & 18437 & 92 & 15235 \\ \hline
151 & 6 & 3800 & 18,24,32,26,25,26 &  & 2 & 93 & 15200 & 15453 & 93 & 15453
\\ \hline
\end{array}
$

\newpage
\end{center}

\subsection{Tables for Section \ref{sec_ParamCyclicConfig}}

\begin{center}
\textbf{Table 6.1. }Values of $v$ for which a cyclic symmetric configuration
$v_{k}$ exists, $5\leq k\leq 51,$ $P(k)\leq v\leq G(k)-1\smallskip $

$\renewcommand{\arraystretch}{1.0}
\begin{array}{@{}r|r|c|c|r|r|r@{}}
\hline
k~ & P(k) & v_{\delta }(k) & v_{\delta }(k)\leq v\leq G(k)-1
\vphantom{L^{L^{L^{L}}}} & E_{c}(k) & G(k) & \text{filling} \\ \hline
5\centerdot & 21 & 21\centerdot & 21_{a},\overline{\emph{22}}_{t}
\vphantom{L^{L^{L^{L}}}} & 23\centerdot & 23 & 100\% \\ \hline
6\centerdot & 31 & 31\centerdot & 31_{a},\overline{\emph{32}}_{t},\overline{
\emph{33}},\overline{\emph{34}^{c}}\vphantom{L^{L^{L^{L}}}} & 35\centerdot &
35 & 100\% \\ \hline
7\centerdot & 43 & 48\centerdot & 48_{b},49,50\vphantom{L^{L^{L}}} &
48\centerdot & 51 & 100\% \\ \hline
8\centerdot & 57 & 57\centerdot & 57_{a},\overline{\emph{58}}_{t},\overline{
\emph{59}^{c}}-\overline{\emph{62}^{c}},63_{b},64-68\vphantom{L^{L^{L^{L}}}}
& 63\centerdot & 69 & 100\% \\ \hline
9\centerdot & 73 & 73\centerdot & 73_{a},\overline{\emph{74}}_{t},\overline{
\emph{75}^{c}}-\overline{\emph{79}^{c}},80_{b},\overline{\emph{81}^{c}}-
\overline{\emph{84}^{c}},85-88\vphantom{L^{L^{L^{L}}}} & 85\centerdot & 89 &
100\% \\ \hline
10\centerdot & 91 & 91\centerdot & 91_{a},\overline{\emph{92}}_{t},\overline{
\emph{93}^{c}}-\overline{\emph{106}^{c}},107-109,110_{c}
\vphantom{L^{L^{L^{L}}}}\smallskip & 107\centerdot & 111 & 100\% \\ \hline
11\centerdot & 111 & 120\centerdot & 120_{b},\overline{\emph{121}^{c}}-
\overline{\emph{132}^{c}},133_{a},\overline{\emph{134}^{c}},135-144
\vphantom{L^{L^{L^{L}}}}\smallskip & 135\centerdot & 145 & 100\% \\ \hline
12\centerdot & 133 & 133\centerdot &
\begin{array}{c}
133_{a},\overline{\emph{134}}_{t},\overline{\emph{135}^{c}}-\overline{\emph{
\ 155}^{c}},156_{c},\overline{\emph{157}^{c}},158,\vphantom{L^{L^{L^{L}}}}
\smallskip \\
159,\overline{\emph{160}^{c}},168_{b},161-170\vphantom{L^{L^{L^{L}}}}
\smallskip
\end{array}
& 161\centerdot & 171 & 100\% \\ \hline
13\centerdot & 157 & 168\centerdot &
\begin{array}{c}
168_{b},\overline{\emph{169}^{c}}-\overline{\emph{182}^{c}},183_{a},
\overline{\emph{184}^{c}}-\overline{\emph{192}^{c}},\vphantom{L^{L^{L^{L}}}}
\smallskip \\
193-212
\end{array}
& 193\centerdot & 213 & 100\% \\ \hline
14\centerdot & 183 & 183\centerdot & 183_{a},\overline{\emph{184}}_{t},
\overline{\emph{185}^{c}}-\overline{\emph{224}^{c}},225-254
\vphantom{L^{L^{L^{L}}}}\smallskip & 225\centerdot & 255 & 100\% \\ \hline
15 & 211 & 255\centerdot &
\begin{array}{c}
255_{b},\overline{\emph{256}^{c}}-\overline{\emph{260}^{c}},\overline{\emph{263}^{c}},272_{c},273_{a},288_{b},\vphantom{L^{L^{L^{L}}}}\smallskip \\
267-302
\end{array}
& 267 & 303 & 95\% \\ \hline
16 & 241 & 255\centerdot &
\begin{array}{c}
255_{b},272_{c},273_{a},288_{b},307_{a},\mathbf{313,317,318,}
\vphantom{L^{L^{L}}} \\
\mathbf{320}-\mathbf{354}
\end{array}
& \mathbf{320} & 355 & 43\% \\ \hline
17 & 273 & 273\centerdot &
\begin{array}{c}
273_{a},\overline{\emph{274}}_{t},288_{b},307_{a},342_{c},\mathbf{343,349,353,}
\vphantom{L^{L^{L^{L}}}} \\
360_{b},381_{a},\mathbf{356}-\mathbf{398}
\end{array}
& \mathbf{356} & 399 & 42\% \\ \hline
18 & 307 & 307\centerdot &
\begin{array}{c}
307_{a},342_{c},360_{b},381_{a},\mathbf{389,391,}
\vphantom{L^{L^{L^{L}}}} \\
\mathbf{395-398,401,}\mathbf{403}-\mathbf{432}
\end{array}
& \mathbf{403} & 433 & 32\% \\ \hline
19 & 343 & \geq 345_{s,t} &
\begin{array}{c}
360_{b},381_{a},\mathbf{445,450,453,455}-\mathbf{458,}\vphantom{L^{L^{L}}} \\
\mathbf{460}-\mathbf{492}\end{array}
& \mathbf{460} & 493 & 29\% \\ \hline
20 & 381 & 381\centerdot &
\begin{array}{c}
381_{a},\overline{\emph{382}},\mathbf{482,497,498,501-503},506_{c},\vphantom{L^{L^{L^{L}}}} \\
\mathbf{505-509,}528_{b},553_{a},\mathbf{511-566}
\vphantom{L^{L^{L^{L}}}}
\end{array}
& \mathbf{511} & 567 & 37\% \\ \hline
21 & 421 & \geq 423_{s,t} &
\begin{array}{c}
506_{c},528_{b},553_{a},\mathbf{586,589,591,592,594,}\vphantom{L^{L^{L}}} \\
\mathbf{595,597,598,}624_{b},\mathbf{600-666}\vphantom{L^{L^{L}}}
\end{array}
& \mathbf{600} & 667 & 32\% \\ \hline
22 & 463 & \geq 465_{br,t} &
\begin{array}{c}
506_{c},528_{b},553_{a},624_{b},\mathbf{633,637},\vphantom{L^{L^{L}}}
\\
\mathbf{640-642},651_{a},\mathbf{644-712}
\end{array}
& \mathbf{644} & 713 & 33\% \\ \hline
\end{array}
$
\end{center}

\noindent Key to Table 6.1: $a\rightarrow
(\ref{eq2_cyclicPG(2,q)})$, $ b\rightarrow
(\ref{eq2_cyclicAG(2,q)}),$ $c\rightarrow
(\ref{eq2_cyclicRuzsa} ),$ $br\rightarrow $Bruck-Ryser Theorem,
$s\rightarrow $\cite{BaumertGordon}, $t\rightarrow $Theorem
\ref{th5_deficiency1}$\vphantom{L^{L}}$\newline
\newpage

\begin{center}
\textbf{Table 6.1} (continue 1) Values of $v$ for which a
cyclic symmetric configuration $v_{k}$ exists, $5\leq k\leq
51,$ $P(k)\leq v<G(k)\smallskip $

$\renewcommand{\arraystretch}{0.91}
\begin{array}{@{}r|c|c|c|r|r}
\hline
k~ & P(k) & v_{\delta }(k) & v_{\delta }(k)\leq v\leq G(k)-1
\vphantom{L^{L^{L^{L}}}} & E_{c}(k) & G(k) \\ \hline
23 & 507 & \geq 509_{br,t} &
\begin{array}{c}
528_{b},553_{a},624_{b},651_{a},
\mathbf{683,686-688,692,}\vphantom{L^{L^{L^{L}}}} \\
 \mathbf{695-700,}728_{b},\mathbf{702-744}\vphantom{L^{L^{L^{L}}}}
\end{array}
& \mathbf{702} & 745 \\ \hline
24 & 553 & 553\centerdot &
\begin{array}{c}
553_{a},\overline{\emph{554}}_{t},624_{b},651_{a},728_{b},\mathbf{738,739,742,}\vphantom{L^{L^{L^{L}}}} \\
\mathbf{747-749,752,753,755,}757_{a}, \vphantom{L^{L^{L^{L}}}} \\
812_{c},840_{b},\mathbf{757}-\mathbf{850}\vphantom{L^{L^{L^{L}}}}
\end{array}
& \mathbf{757} & 851 \\ \hline
25 & 601 & \geq 602_{s} &
\begin{array}{c}
624_{b},651_{a},728_{b},757_{a},812_{c},\mathbf{830},840_{b},
\vphantom{L^{L^{L}}} \\
871_{a},930_{c}, 960_{b},\mathbf{837}-\mathbf{960} \vphantom{L^{L^{L}}}
\end{array}
& \mathbf{837} & 961 \\ \hline
26 & 651 & 651\centerdot &
\begin{array}{c}
651_{a},\overline{\emph{652}}_{t},728_{b},757_{a},812_{c},840_{b},871_{a},
\vphantom{L^{L^{L^{L}}}} \\
\mathbf{885,888,895,900,903,905-907},\vphantom{L^{L^{L}}}\\
\mathbf{910-913,915-917,919-925,} \vphantom{L^{L^{L}}}\\
\mathbf{927,}930_{c},960_{b}, \mathbf{929-984}
\end{array}
& \mathbf{929} & 985 \\ \hline
27 & 703 & \geq 704_{s} &
\begin{array}{c}
728_{b},757_{a},812_{c},840_{b},871_{a},930_{c},960_{b},\vphantom{L^{L^{L}}}
\\
\mathbf{970,971,972,975,977},\mathbf{978,}\mathbf{985,987},\mathbf{988,} \vphantom{L^{L^{L}}}\\
\mathbf{991},993_{a},\mathbf{993-997,1000,1001,} \\
\mathbf{1003}-\mathbf{1015,}1023_{b},1057_{a},\mathbf{1017}-\mathbf{1106}
\end{array}
& \mathbf{1017} & 1107 \\ \hline
28 & 757 & 757\centerdot &
\begin{array}{c}
757_{a},\overline{\emph{758}}_{t},812_{c},840_{b},871_{a},930_{c},960_{b},
993_{a},\vphantom{L^{L^{L^{L}}}} \\
\mathbf{1006},1023_{b},\mathbf{1045,1051,1053,}1057_{a}, \\
\mathbf{1063-1067,1070-1972,}\mathbf{1074,1075,} \\
\mathbf{1077,1079}-\mathbf{1170}
\end{array}
& \mathbf{1079} & 1171 \\ \hline
29 & 813 & \geq 815_{s,t} &
\begin{array}{c}
840_{b},871_{a},930_{c},960_{b},993_{a},1023_{b},1057_{a},
\vphantom{L^{L^{L^{L}}}} \\
\mathbf{1091,1127,1135,1137,1141,1143,} \vphantom{L^{L^{L^{L}}}}\\
\mathbf{1145,1146,1151}-\mathbf{1246}
\end{array}
& \mathbf{1151} & 1247 \\ \hline
30 & 871 & 871\centerdot &
\begin{array}{c}
871_{a},\overline{\emph{872}}_{t},930_{c},960_{b},993_{a},1023_{b},1057_{a},
\vphantom{L^{L^{L^{L}}}} \\
\mathbf{1196,1198-1201,1206,1207,1216,} \vphantom{L^{L^{L^{L}}}}\\
\mathbf{1217,1219-1224,}1332_{c},\mathbf{1226-1360}
\end{array}
& \mathbf{1226} & 1361 \\ \hline
31 & 931 & \geq 933_{br,t} &
\begin{array}{c}
960_{b},993_{a},1023_{b},1057_{a},\mathbf{1298,1309,}\vphantom{L^{L^{L^{L}}}}
\\
\mathbf{1314,1315,1320,1321,1324,1325,} \\
1332_{c},\mathbf{1330-1335,}\\
\mathbf{1339-1346,}1368_{b},\mathbf{1348-1494}\end{array}
& \mathbf{1348} & 1495 \\ \hline
32 & 993 & 993\centerdot &
\begin{array}{c}
993_{a},\overline{\emph{994}}_{t},1023_{b},1057_{a},1332_{c},\mathbf{1366,}1368_{b},
\vphantom{L^{L^{L^{L}}}} \\
\mathbf{1383,1388,1391-1395,1397,1398,} \\
\mathbf{1400,1401,1403,}1407_{a},\\
\mathbf{1406-1409,}\mathbf{1411-1414,1416,} \\
\mathbf{1420,1421,1424-1434,1436-1568}
\end{array}
& \mathbf{1436} & 1569 \\ \hline
\end{array}
$
\end{center}

\noindent Key to Table 6.1: $a\rightarrow
(\ref{eq2_cyclicPG(2,q)})$, $ b\rightarrow
(\ref{eq2_cyclicAG(2,q)}),$ $c\rightarrow
(\ref{eq2_cyclicRuzsa} ),\vphantom{L^{L}}$ $br\rightarrow
$Bruck-Ryser Theorem, $s\rightarrow $\cite{BaumertGordon},
$t\rightarrow $Theorem \ref{th5_deficiency1}\newpage

\begin{center}
\textbf{Table 6.1} (continue 2) Values of $v$ for which a
cyclic symmetric configuration $v_{k}$ exists, $5\leq k\leq
51,$ $P(k)\leq v<G(k)\smallskip $

$
\begin{array}{@{}r|c|c|c|r|r|}
\hline
k~ & P(k) & v_{\delta }(k) & v_{\delta }(k)\leq v\leq G(k)-1
\vphantom{L^{L^{L^{L}}}} & E_{c}(k) & G(k) \\ \hline
33 & 1057 & 1057\centerdot &
\begin{array}{c}
1057_{a},\overline{\emph{1058}}_{t},1332_{c},1368_{b},1407_{a},\mathbf{\
1492,1506,}\vphantom{L^{L^{L^{L}}}} \\
\mathbf{1507,1515,1518,1520,1521,1528,1529,} \\
\mathbf{1533,1535,1537,1540,1542,1543,1545,} \\
\mathbf{1547-1553,}\mathbf{1555}-\mathbf{1559},\\
1640_{c},1680_{b},\mathbf{1561}-\mathbf{1718}
\end{array}
& \mathbf{1561} & 1719 \\ \hline
34 & 1123 & \geq 1125_{br,t} &
\begin{array}{c}
1332_{c},1368_{b},1407_{a},1640_{c},\mathbf{1664,1665,1670,}
\vphantom{L^{L^{L^{L}}}} \\
\mathbf{1676},1680_{b},\mathbf{1686,1693},\mathbf{1698,1699,1702,} \\
\mathbf{1705,1708}-\mathbf{1712,1714,1717,1721,} \\
1723_{a},\mathbf{1723}-\mathbf{1726,1728,1730}-\mathbf{1742,} \\
\mathbf{1744}-\mathbf{1752,}1806_{c},1848_{b},\mathbf{1754}-\mathbf{1876}
\end{array}
& \mathbf{1754} & 1877 \\ \hline
35 & 1191 & \geq 1193_{s,t} &
\begin{array}{c}
1332_{c},1368_{b},1407_{a},1640_{c},1680_{b},1723_{a},
\vphantom{L^{L^{L^{L}}}} \\
\mathbf{1777,1781,1783,1788,1792,1793,1795,} \\
\mathbf{1798,1800}-\mathbf{1803,}1806_{c},\mathbf{1805}-\mathbf{1807,} \\
\mathbf{1810,1812}-\mathbf{1815,}1848_{b},1893_{a}, \\
\mathbf{1817}-\mathbf{1974}
\end{array}
& \mathbf{1817} & 1975 \\ \hline
36 & 1261 & \geq 1262_{s} &
\begin{array}{c}
1332_{c},1368_{b},1407_{a},1640_{c},1680_{b},1723_{a},1806_{c},
\vphantom{L^{L^{L^{L}}}} \\
1848_{b},\mathbf{1853,1855,1860,1867-1870,}\\
\mathbf{1872-1876,1878,1882-1884,} 1893_{a},\\
\mathbf{1886-2010}
\end{array}
& \mathbf{1886} & 2011 \\ \hline
37 & 1333 & \geq 1335_{s,t} &
\begin{array}{c}
1368_{b},1407_{a},1640_{c},1680_{b},1723_{a},1806_{c},1848_{b},
\vphantom{L^{L^{L^{L}}}} \\
\mathbf{1892,}1893_{a},\mathbf{1910,1922,1930,1934,1938,}\\
\mathbf{1943,1944,1947-1953,1957,1959,} \\
\mathbf{1960,1962,1963,1965-1967,1969,} \\
2162_{c},\mathbf{1972-2198}
\end{array}
& \mathbf{1972} & 2199 \\ \hline
38 & 1407 & 1407\centerdot &
\begin{array}{c}
1407_{a},1640_{c},1680_{b},1723_{a},1806_{c},1848_{b},
\vphantom{L^{L^{L^{L}}}} \\
1893_{a},\mathbf{2059,2061,2073,2088,2089,} \\
\mathbf{2092,2094,2096,2097,2099-2101,} \\
\mathbf{2103,2105-2108,2110,2111,} \\
\mathbf{2114}-\mathbf{2116,2118,2123,2124,} \\
\mathbf{2126}-\mathbf{2130,2135}-\mathbf{2153,2155}\text{\textbf{-}}\mathbf{
\ 2157,} \\
2162_{c},\mathbf{2159-2164},\mathbf{2166-2170}, \\
2208_{b},2257_{a}\mathbf{,2172-2292}
\end{array}
& \mathbf{2172} & 2293 \\ \hline
\end{array}
$
\end{center}

\noindent Key to Table 6.1: $a\rightarrow
(\ref{eq2_cyclicPG(2,q)})$, $ b\rightarrow
(\ref{eq2_cyclicAG(2,q)}),$ $c\rightarrow
(\ref{eq2_cyclicRuzsa} ),\vphantom{L^{L}}$ $br\rightarrow
$Bruck-Ryser Theorem, $s\rightarrow $\cite{BaumertGordon},
$t\rightarrow $Theorem \ref{th5_deficiency1}\newpage

\begin{center}
\textbf{Table 6.1} (continue 3) Values of $v$ for which a
cyclic symmetric configuration $v_{k}$ exists, $5\leq k\leq
51,$ $P(k)\leq v<G(k)\smallskip $

$
\begin{array}{@{}r|c|c|c|r|r|}
\hline
k~ & P(k) & v_{\delta }(k) & v_{\delta }(k)\leq v\leq G(k)-1
\vphantom{L^{L^{L^{L}}}} & E_{c}(k) & G(k) \\ \hline
39 & 1483 & \geq 1485_{br,t} &
\begin{array}{c}
1640_{c},1680_{b},1723_{a},1806_{c},1848_{b},1893_{a},
\vphantom{L^{L^{L^{L}}}} \\
2162_{c},\mathbf{2187,2195,}2208_{b},\mathbf{2240,2241,} \\
\mathbf{2243,2247,2248,2251,2252,2254,} \\
\mathbf{2255,}2257_{a},\mathbf{2258-2260,2263-2265,} \\
\mathbf{2269,2270,2277-2279,2281,2283,} \\
\mathbf{2287,2289-2309,2311,2313-2328,} \\
2400_{b}\mathbf{,}2451_{a}\mathbf{,2330-2504}
\end{array}
& \mathbf{2330} & 2505 \\ \hline
40 & 1561 & \geq 1563_{s,t} &
\begin{array}{c}
1640_{c},1680_{b},1723_{a},1806_{c},1848_{b},1893_{a},
\vphantom{L^{L^{L^{L}}}} \\
2162_{c},2208_{b},2257_{a},\mathbf{2326,2338,2345,} \\
\mathbf{2349,2353,2355,2357,2360,2361,} \\
\mathbf{2363,2364,2366-2368,2370,2372,} \\
\mathbf{2374-2377,2379-2381,2387-2389,} \\
\mathbf{2393,2395-2397,2399,}2400_{b}\mathbf{,2401,} \\
\mathbf{2403,2404,2406,2407,2409,} \\
\mathbf{2411-2418,}2451_{a}\mathbf{,2420-2436,} \\
\mathbf{2438-2564}
\end{array}
& \mathbf{2438} & 2565 \\ \hline
41 & 1641 & \geq 1643_{s,t} &
\begin{array}{c}
1680_{b},1723_{a},1806_{c},1848_{b},1893_{a},2162_{c},
\vphantom{L^{L^{L^{L}}}} \\
2208_{b},2257_{a},\mathbf{2345,2399,}2400_{b},\mathbf{2436,} \\
\mathbf{2449,}2451_{a}\mathbf{,2459,2460,2465,2471,} \\
\mathbf{2472,2479,2480,2483,2485,2491,} \\
\mathbf{2493,2494,2496-2500,2502,2503,} \\
\mathbf{2505,2507-2513,2515-2525,} \\
\mathbf{2528-2540,2542,2544}-\mathbf{2610}
\end{array}
& \mathbf{2544} & 2611 \\ \hline
42 & 1723 & 1723\centerdot &
\begin{array}{c}
1723_{a},\overline{\emph{1724}}_{t}\emph{,}
1806_{c},1848_{b},1893_{a},2162_{c},\vphantom{L^{L^{L^{L}}}} \\
2208_{b},2257_{a},2400_{b},2451_{a}\mathbf{,2510,2522,} \\
\mathbf{2539,2541,2557-2559,2562,2564,} \\
\mathbf{2566-2568,2570,2573,2577,2578,} \\
\mathbf{2580-2584,2586-2590,2593,2595,} \\
\mathbf{2597-2601,2603-2610,2612,} \\
\mathbf{2613,2615-2626,}2756_{c},\mathbf{2628}-\mathbf{2794}
\end{array}
& \mathbf{2628} & 2795 \\ \hline
\end{array}
$
\end{center}

\noindent Key to Table 6.1: $a\rightarrow
(\ref{eq2_cyclicPG(2,q)})$, $ b\rightarrow
(\ref{eq2_cyclicAG(2,q)}),$ $c\rightarrow
(\ref{eq2_cyclicRuzsa} ),\vphantom{L^{L}}$ $br\rightarrow
$Bruck-Ryser Theorem, $s\rightarrow $\cite{BaumertGordon},
$t\rightarrow $Theorem \ref{th5_deficiency1}\newpage

\begin{center}
\textbf{Table 6.1} (continue 4) Values of $v$ for which a
cyclic symmetric configuration $v_{k}$ exists, $5\leq k\leq
51,$ $P(k)\leq v<G(k)\smallskip $

$
\begin{array}{@{}r|c|c|c|r|r|}
\hline
k~ & P(k) & v_{\delta }(k) & v_{\delta }(k)\leq v\leq G(k)-1
\vphantom{L^{L^{L^{L}}}} & E_{c}(k) & G(k) \\ \hline
43 & 1807 & \geq 1809_{br,t} &
\begin{array}{c}
1848_{b},1893_{a},2162_{c},2208_{b},2257_{a},2400_{b},
\vphantom{L^{L^{L^{L}}}} \\
2451_{a}\mathbf{,2684,2686,2688,2715,2725,} \\
\mathbf{2728,2734,2737,2739,2744,2752,} \\
2756_{c},\mathbf{2757,2759,2762,2763,} \\
\mathbf{2766-2768,2771,2772,2776,2777,} \\
\mathbf{2783,2786-2789,2791,2792,} \\
\mathbf{2794-2798,2800,2801,}2808_{b}, \\
\mathbf{2803-2811,2813-2815,2817-2858,} \\
2863_{a},\mathbf{2860-3014}
\end{array}
& \mathbf{2860} & 3015 \\ \hline
44 & 1893 & 1893\centerdot &
\begin{array}{c}
1893_{a},\overline{\emph{1894}}_{t}\emph{,}
2162_{c},2208_{b},2257_{a},2400_{b},\vphantom{L^{L^{L^{L}}}} \\
2451_{a}\mathbf{,}2756_{c},2808_{b},\mathbf{2811,2821,2826,} \\
\mathbf{2834,2836,2844,2848,2849,2861,} \\
\mathbf{2862,}2863_{a},\mathbf{2865-2867,2870,2871,} \\
\mathbf{2873-2875,2879-2881,2884,2887,} \\
\mathbf{2890-2895,2898,2899,2901-2912,} \\
\mathbf{2914,2916-3192}
\end{array}
& \mathbf{2916} & 3193 \\ \hline
45 & 1981 & \geq 1983_{s,t} &
\begin{array}{c}
2162_{c},2208_{b},2257_{a},2400_{b},2451_{a},2756_{c},
\vphantom{L^{L^{L^{L}}}} \\
2808_{b},2863_{a},\mathbf{2994,3013,3014,3019,} \\
\mathbf{3038,3052,3054,3056,3066,3069,} \\
\mathbf{3082-3085,3087,3088,3090,3091,} \\
\mathbf{3093,3095-3098,3101-3103,} \\
\mathbf{3105-3112,3114,3116-3120,} \\
\mathbf{3122-3163,3165-3374}
\end{array}
& \mathbf{3165} & 3375 \\ \hline
46 & 2071 & \geq 2073_{s,t} &
\begin{array}{c}
2162_{c},2208_{b},2257_{a},2400_{b},2451_{a},2756_{c},
\vphantom{L^{L^{L^{L}}}} \\
2808_{b},2863_{a},\mathbf{3124,3171,3188,3191,} \\
\mathbf{3194,3196,3197,3198,3206,3216,} \\
\mathbf{3218,3219,3221,3227,3228,3231,} \\
\mathbf{3233,3234,3238-3241,3244-3247,} \\
\mathbf{3249,3250,3252-3265,3267-3271,} \\
\mathbf{3273-3278,3280-3406}
\end{array}
& \mathbf{3280} & 3407 \\ \hline
\end{array}
$
\end{center}

\noindent Key to Table 6.1: $a\rightarrow
(\ref{eq2_cyclicPG(2,q)})$, $ b\rightarrow
(\ref{eq2_cyclicAG(2,q)}),$ $c\rightarrow
(\ref{eq2_cyclicRuzsa} ),\vphantom{L^{L}}$ $br\rightarrow
$Bruck-Ryser Theorem, $s\rightarrow $\cite{BaumertGordon},
$t\rightarrow $Theorem \ref{th5_deficiency1}\newpage

\begin{center}
\textbf{Table 6.1 (}continue 5)\textbf{\ }Values of $v$ for
which a cyclic symmetric configuration $v_{k}$ exists, $5\leq
k\leq 51,$ $P(k)\leq v<G(k)\smallskip $

$
\begin{array}{@{}r|c|c|c|r|r|}
\hline
k~ & P(k) & v_{\delta }(k) & v_{\delta }(k)\leq v\leq G(k)-1
\vphantom{L^{L^{L^{L}}}} & E_{c}(k) & G(k) \\ \hline
47 & 2163 & 2165_{br,t} &
\begin{array}{c}
2208_{b},2257_{a},2400_{b},2451_{a},2756_{c},2808_{b},
\vphantom{L^{L^{L^{L}}}} \\
2863_{a},\mathbf{3255,3261,3271,3280,3285,} \\
\mathbf{3292,3301,3312,3327,3331,3342,} \\
\mathbf{3343,3346-3348,3351,3353,} \\
\mathbf{3355-3357,3360-3371,3375-3379,} \\
\mathbf{3381-3384,3387,3389,3390,} \\
\mathbf{3392-3402,3404-3408,3412,3413,} \\
3422_{c},\mathbf{3415-3427,}3480_{b},3541_{a}\mathbf{,} \\
\mathbf{3429-3608}
\end{array}
& \mathbf{3429} & 3609 \\ \hline
48 & 2257 & 2257\centerdot &
\begin{array}{c}
2257_{a},\overline{\emph{2258}}_{t}\emph{,}
2400_{b},2451_{a},2756_{c},2808_{b},\vphantom{L^{L^{L^{L}}}} \\
2863_{a},\mathbf{3418,}3422_{c},\mathbf{3431,3445,3459,} \\
3480_{b},\mathbf{3487,3491,3492,3495,3499,} \\
\mathbf{3509,3512,3515,3518,3519,3522,} \\
\mathbf{3523,3526-3528,3535,3540,} \\
3541_{a},\mathbf{3545-3547,3549-3552,} \\
\mathbf{3554-3563,3565,3567,3568,3570,} \\
\mathbf{3572-3578,3580-3583,3586,3587,} \\
\mathbf{3589-3595,3597-3605,3607-3618,} \\
\mathbf{3620-3630,3632-3640,}3660_{c}, \\
\mathbf{3642-3774}
\end{array}
& \mathbf{3642} & 3775 \\ \hline
49 & 2353 & \geq 2354_{s} &
\begin{array}{c}
2400_{b},2451_{a},2756_{c},2808_{b},2863_{a},3422_{c},
\vphantom{L^{L^{L^{L}}}} \\
3480_{b},3541_{a},\mathbf{3608,3627,3637,3640,} \\
\mathbf{3642,3644,3646,3647,3649,3652,} \\
\mathbf{3653,3655,}3660_{c},\mathbf{3665,3669,3671,} \\
\mathbf{3675,3677,3678,3680-3683,3685,} \\
\mathbf{3688-3692,3694-3701,3707,3709,} \\
3720_{b},\mathbf{3711-3722,3724-3730,} \\
\mathbf{3732-3735,3737-3740,3742,3743,} \\
\mathbf{3745,3747-3755,}3783_{a},\mathbf{3757-3825,} \\
\mathbf{3827-3837,3839-3916}
\end{array}
& \mathbf{3839} & 3917 \\ \hline
\end{array}
$
\end{center}

\noindent Key to Table 6.1: $a\rightarrow
(\ref{eq2_cyclicPG(2,q)})$, $ b\rightarrow
(\ref{eq2_cyclicAG(2,q)}),$ $c\rightarrow
(\ref{eq2_cyclicRuzsa} ),\vphantom{L^{L}}$ $br\rightarrow
$Bruck-Ryser Theorem, $s\rightarrow $\cite{BaumertGordon},
$t\rightarrow $Theorem \ref{th5_deficiency1}\newpage

\begin{center}
\textbf{Table 6.1} (continue 6). Values of $v$ for which a cyclic symmetric
configuration $v_{k}$ exists, $5\leq k\leq 51,$ $P(k)\leq v<G(k)\smallskip $

$
\begin{array}{@{}r|c|c|c|r|r|}
\hline
k~ & P(k) & v_{\delta }(k) & v_{\delta }(k)\leq v\leq G(k)-1
\vphantom{L^{L^{L^{L}}}} & E_{c}(k) & G(k) \\ \hline
50 & 2451 & 2451\centerdot &
\begin{array}{c}
2451_{a},\overline{\emph{2452}}_{t}\emph{,}
2756_{c},2808_{b},2863_{a},3422_{c},\vphantom{L^{L^{L^{L}}}} \\
3480_{b},3541_{a},3660_{c},\mathbf{3685,3688,3692,} \\
\mathbf{3712,3714,3716,}3720_{b},\mathbf{3726,3743,} \\
\mathbf{3745,3749,3750,3752,3753,3758,} \\
\mathbf{3762,3766,3767,3769,3770,3772,} \\
\mathbf{3775,3779,3780,}3783_{a},\mathbf{3782-3785,} \\
\mathbf{3788-3791,3796,3803,3805,3811,} \\
\mathbf{3813,3817-3823,3825,} \\
\mathbf{3828-3832,3834-3840,3842,3843,} \\
\mathbf{3845,3846,3848-3869,}4095_{b},4161_{a}, \\
\mathbf{3871-4188}
\end{array}
& \mathbf{3871} & 4189 \\ \hline
51 & 2551 & \geq 2552_{s} &
\begin{array}{c}
2756_{c},2808_{b},2863_{a},3422_{c},3480_{b},3541_{a},
\vphantom{L^{L^{L^{L}}}} \\
3660_{c},3720_{b},3783_{a},\mathbf{3871,3894,3927,} \\
\mathbf{3938,3986,3998,4004,4018,4032,} \\
\mathbf{4035,4037,4042,4048-4050,4053,} \\
\mathbf{4060,4064-4066,4068,4072,4076,} \\
\mathbf{4079-4083,4085,4087,4090,4091,} \\
\mathbf{4093,}4095_{b},\mathbf{4095-4100,4102-4104,} \\
\mathbf{4107-4109,4112-4120,4124,4125,} \\
\mathbf{4127-4131,4133-4156,}4161_{a}, \\
\mathbf{4158-4164,4166-4196,4199,4200,} \\
\mathbf{4202-4206,4208-4231,4233-4380}
\end{array}
& \mathbf{4233} & 4381 \\ \hline
\end{array}
$
\end{center}

\noindent Key to Table 6.1: $a\rightarrow
(\ref{eq2_cyclicPG(2,q)})$, $ b\rightarrow
(\ref{eq2_cyclicAG(2,q)}),$ $c\rightarrow
(\ref{eq2_cyclicRuzsa} ),\vphantom{L^{L}}$ $br\rightarrow
$Bruck-Ryser Theorem, $s\rightarrow $\cite{BaumertGordon},
$t\rightarrow $Theorem \ref{th5_deficiency1}

\begin{center}
\textbf{Table 6.2}. Upper bounds on the cyclic existence bound
$E_{c}(k),$ $ 52\leq k\leq 83\bigskip $

$\renewcommand{\arraystretch}{1.0}
\begin{array}{|ccc|ccc|crr|ccc|}
\hline
k & E_{c}(k) & G(k) & k & E_{c}(k) & G(k) & k & E_{c}(k) & G(k) & k &
E_{c}(k) & G(k) \\ \hline
52 & \mathbf{4359} & 4541 & 60 & \mathbf{5687} & 6039 & 68 & \mathbf{7463} &
7913 & 76 & \mathbf{10023} & 10179 \\
53 & \mathbf{4463} & 4695 & 61 & \mathbf{5994} & 6269 & 69 & \mathbf{8111} &
8291 & 77 & \mathbf{10229} & 10409 \\
54 & \mathbf{4513} & 4747 & 62 & \mathbf{6150} & 6431 & 70 & \mathbf{8125} &
8435 & 78 & \mathbf{10395} & 10599 \\
55 & \mathbf{5195} & 5197 & 63 & \mathbf{6611} & 6783 & 71 & \mathbf{8288} &
8661 & 79 & \mathbf{10800} & 10817 \\
56 & \mathbf{5341} & 5451 & 64 & \mathbf{6796} & 7055 & 72 & \mathbf{8694} &
8947 & 80 & \mathbf{10977} & 11127 \\
57 & \mathbf{5501} & 5547 & 65 & \mathbf{6853} & 7187 & 73 & \mathbf{8813} &
9027 & 81 & \mathbf{11396} & 11435 \\
58 & \mathbf{5551} & 5703 & 66 & \mathbf{7279} & 7515 & 74 & \mathbf{8965} &
9507 & 82 & \mathbf{11443} & 11629 \\
59 & \mathbf{5612} & 5823 & 67 & \mathbf{7359} & 7639 & 75 & \mathbf{9883} &
9965 & 83 & \mathbf{11593} & 12041 \\ \hline
\end{array}
$\newpage
\end{center}

\begin{center}
\textbf{Table 6.3}. New cyclic symmetric configurations $
v_{k}$
\begin{tabular}{@{}c|@{\,\,~}c|@{\,\,~}c@{}}
\hline
  $k$&$v$&the 1-st row of incidence matrix  \\    \hline
  17 & 382& 0,25,69,81,88,89,112,123,126,128,141,174,196,202,206,223,232 \\
  17 &$
  \begin{array}{c}
  383\\
   385\\
   386\\
   388
   \end{array}$
   & 0,5,15,34,35,42,73,75,86,89,98,134,151,155,177,183,201 \\
  17 & 384 &0,68,70,84,90,107,111,120,139,151,185,186,193,196,211,244,249 \\
  17 & 387 & 0,5,7,17,52,56,67,80,81,100,122,138,159,165,168,191,199 \\
  18&389&0,10,27,28,35,50,74,94,103,105,108,146,159,165,191,195,207,228 \\
  18&391&0, 1, 4,15,35,42,75,85,94,111,133,139,141,157,162,194,206,219 \\
  18&395&0,71,73,81,85,88 117 118 141 167 172 183 192 210 231 244 250 272\\
    18&396&0,79,89,106,107,114,129,153,173,182,184,187,225,238,244,270,274,286 \\
 18&397&0,36,50,56,59,74,78,85,122,139,147,149,179,180,192,226,231,247 \\
 18&398&0,18,30,32,71,84,90,93,119,127,152,169,176,192,196,197,207,243 \\
  18& 401 & 0,71,94,104,136,160,164,176,195,217,238,243,256,263,290,292,293,301 \\
  18&404&0, 2,10,14,17,46,47,70,96,101,112,121,139,160,173,179,201,236\\
  18&$
  \begin{array}{c}
  405\\
    407
   \end{array}$    & 0,2,10,22,53,56,82,83,89,98,130,148,153,167,188,192,205,216 \\
  18& 406 & 0,49,59,62,97,99,117,141,173,180,184,192,201,206,207,237,253,278 \\
  18&408&0,14,15,27,50,53,60,81,97,115,137,139,145,156,188,208,213,217\\
  18& 409 & 0,91,193,195,203,215,246,249,275,276,282,291,323,341,346,360,381,385 \\
18  &410&0,36,68,103,110,111,121,124,130,161,176,200,202,225,230,247,259,263\\
  18&411&0,93,195,197,205,217,248,251,277,278,284,293,325,343,348,362,383,387 \\
  \hline
\end{tabular}
\end{center}
\newpage
\subsection{Tables for Section \ref{sec6_spectrum}}

\begin{center}
\textbf{Table 7.1.} The currently known parameters of symmetric configurations
$v_{k}$ (cyclic and non-cyclic)\smallskip

$\renewcommand{\arraystretch}{1}
\begin{array}{@{}c|r|c|r|r|r}
\hline
&  &  & E(k) &  &  \\
k & P(k) & P(k)\leq v\leq G(k)-1 & \leq ~ & G(k) & \text{filling} \\ \hline
3\centerdot & 7 & 7\vphantom{L^{L^{L^{L}}}} & 7\centerdot & 7 & 100\% \\
\hline
4\centerdot & 13 & 13\vphantom{L^{L^{L^{L}}}} & 13\centerdot & 13 & 100\% \\
\hline
5\centerdot & 21 & 21_{a},\overline{\emph{22}}_{t}\vphantom{L^{L^{L^{L}}}} &
23\centerdot & 23 & 100\% \\ \hline
6\centerdot & 31 & 31_{a},\overline{\emph{32}}_{t},\overline{\emph{33}},34
\vphantom{L^{L^{L^{L}}}} & 34\centerdot & 35 & 100\% \\ \hline
7 & 43 & \overline{\emph{43}}_{br},\overline{\emph{44}}_{t},45,48_{b\cdot
f\cdot r},49_{m\cdot r},50_{m\cdot r}\vphantom{L^{L^{L^{L}}}} & 48
\phantom{\centerdot} & 51 & 75\% \\ \hline
8 & 57 & 57_{a},\overline{\emph{58}}_{t},63_{b\cdot f\cdot r},
64_{m\cdot r}-68_{m\cdot r}\vphantom{L^{L^{L^{L}}}} & 63\phantom{\centerdot} & 69 & 67\%
\\ \hline
9 & 73 & 73_{a},\overline{\emph{74}}_{t},78_{g\cdot h},80_{b\cdot f\cdot r},
81_{m\cdot r}-88_{m\cdot r}\vphantom{L^{L^{L^{L}}}} & 80
\phantom{\centerdot} & 89 & 75\% \\ \hline
10 & 91 & 91_{a},\overline{\emph{92}}_{t},98_{j},107_{y}-109_{y},110_{c\cdot
f\cdot k\cdot m\cdot r\cdot S\cdot T}\vphantom{L^{L^{L^{L}}}} & 107
\phantom{\centerdot} & 111 & 35\% \\ \hline
11 & 111 &
\begin{array}{c}
\overline{\emph{111}},\overline{\emph{112}},120_{b\cdot f\cdot r},
121_{m\cdot r}-133_{m\cdot r},135_{y}-142_{y},\vphantom{L^{L^{L^{L}}}} \\
143_{m\cdot y\cdot S},144_{m\cdot y\cdot r\cdot S\cdot T}
\vphantom{L^{L^{L^{L}}}}
\end{array}
& 135\phantom{\centerdot} & 145 & 76\% \\ \hline
12 & 133 &
\begin{array}{c}   133_{a},\overline{\emph{134}}_{t},135,154_{\lambda},155_{\lambda},156_{m\cdot r\cdot S\cdot
T}-169_{m\cdot r\cdot S\cdot T},
\vphantom{L^{L^{L^{L}}}} \\
170_{m\cdot r\cdot T\cdot X}\vphantom{L^{L^{L^{L}}}}
\end{array}
 & 154\phantom{\centerdot} & 171 & 52\% \\ \hline
13 & 157 &
\begin{array}{c}
\overline{\emph{158}}_{t},168_{b\cdot f\cdot r},169_{m\cdot r}-183_{m\cdot r},189_{g},193_{y}-209_{y},\vphantom{L^{L^{L^{L}}}} \\
210_{m\cdot y\cdot r}-212_{m\cdot y\cdot r}\vphantom{L^{L^{L^{L}}}}
\end{array}
& 193\phantom{\centerdot} & 213 & 68\% \\ \hline
14 & 183 & 183_{a},\overline{\emph{184}}_{t},210_{g},222_{\lambda},223_{\lambda},224_{m},225_{m\cdot
y\cdot r}-254_{m\cdot y\cdot r}\vphantom{L^{L^{L^{L}}}} & 222
\phantom{\centerdot} & 255 & 50\% \\ \hline
15 & 211 &
\begin{array}{c}\overline{\emph{211}}_{br},\overline{\emph{212}}
_{t},231_{g},238_{\lambda},239_{\lambda},240_{m\cdot r}-266_{m\cdot r}, \vphantom{L^{L^{L^{L}}}}  \\
267_{m\cdot y\cdot r}-302_{m\cdot y\cdot r}\vphantom{L^{L^{L^{L}}}}
\end{array}
 & 238\phantom{\centerdot} & 303 & 73\% \\
\hline
16 & 241 &
\begin{array}{c}
252_{g\cdot h},255_{b\cdot f\cdot r},256_{m\cdot r}-321_{m\cdot r},
322_{\lambda\cdot r\cdot T\cdot Z},\vphantom{L^{L^{L^{L}}}} \\
323_{m\cdot r\cdot S}-354_{m\cdot r\cdot S}\vphantom{L^{L^{L^{L}}}}
\end{array}
&255\phantom{\centerdot} & 355 & 89\% \\ \hline
17 & 273 &
\begin{array}{c}
273_{a},\overline{\emph{274}}_{t},288_{b\cdot f\cdot r},289_{m\cdot r}-307_{m\cdot r},
321_{\lambda},322_{\lambda},\vphantom{L^{L^{L^{L}}}} \\
323_{m\cdot S},372_{W},324_{m\cdot r}-381_{m\cdot r},\mathbf{382}_{\mathbf{Z}}-\mathbf{388}_{
\mathbf{Z}},
\vphantom{L^{L^{L^{L}}}}\\
389_{\lambda\cdot Z},390_{\lambda\cdot Z},391_{m\cdot S\cdot Z}-398_{m\cdot S\cdot Z}\vphantom{L^{L^{L^{L}}}}
\end{array}
& \mathbf{321}\phantom{\centerdot} & 399 & 79\% \\ \hline
18 & 307 &
\begin{array}{c}
307_{a},340_{\lambda},341_{\lambda},342_{m\cdot r}-381_{m\cdot r},
\vphantom{L^{L^{L^{L}}}} \\
\mathbf{389}_{Z},\mathbf{391}_{Z},\mathbf{395}_{\mathbf{Z}}-\mathbf{398}_{\mathbf{Z}},\mathbf{401}_{\mathbf{Z}},403_{g\cdot Z},
\vphantom{L^{L^{L^{L}}}}\\
\mathbf{404}_{\mathbf{Z}}-\mathbf{411}_{\mathbf{Z}},412_{\lambda\cdot Z},413_{\lambda\cdot Z},414_{m\cdot S\cdot Z}-432_{m\cdot S\cdot Z}\vphantom{L^{L^{L^{L}}}}
\end{array}
& \mathbf{403}\phantom{\centerdot} & 433 & 63\% \\ \hline
\end{array}
\medskip $

Key to Table 7.1: $a\rightarrow (\ref{eq2_cyclicPG(2,q)})$,
$b\rightarrow ( \ref{eq2_cyclicAG(2,q)}),$ $c\rightarrow
(\ref{eq2_cyclicRuzsa}),$ $ f\rightarrow
(\ref{eq2_q-1-cancel}),$ $g\rightarrow (\ref{eq2_Baer}),$ $
h\rightarrow (\ref{eq2_Baer2}),$ $j\rightarrow
(\ref{eq2_FuLabNabDecomp}),$ $ k\rightarrow
(\ref{eq2_affine_q-1}),$  $\lambda\rightarrow
\eqref{eq2_Balbuena1}-\eqref{eq2_AaParBalb}$, $m\rightarrow
(\ref{eq2_tetaExten} ),P\rightarrow
(\ref{eq3.3_AfPlnPairsWeitghts}),r\rightarrow (\ref
{eq4_AfPlExten}),S\rightarrow (\ref{eq4_RuzExten})$,
$T\rightarrow (\ref {eq4_RuzExten_d=p-1}),W\rightarrow $
Table~4.1, $X\rightarrow $Table~4.2$ ,y\rightarrow $ Table~6.1
with $k\leq 15$, $Z\rightarrow $ Table~6.1 with $ k>15$,
$br\rightarrow $Bruck-Ryser Theorem, $t\rightarrow $Theorem
\ref{th5_deficiency1} \newpage

\textbf{Table 7.1} (continue 1). The currently known parameters of symmetric
configurations $v_{k}$ (cyclic and non-cyclic)\smallskip

$\renewcommand{\arraystretch}{1}
\begin{array}{@{}r|r|@{}c@{}|@{\,\,}r|r|@{\,\,}r@{}}
\hline
&  &  & E(k) &  &  \\
k & P(k) & P(k)\leq v\leq G(k)-1 & \leq ~ & G(k) & \text{filling} \\ \hline
19 & 343 &
\begin{array}{c}
\overline{\emph{344}}_{t},360_{b\cdot f\cdot r},361_{m\cdot r}-381_{m\cdot r},
434_{g\cdot W},435_{\lambda},436_{\lambda},\vphantom{L^{L^{L^{L}}}} \\
437_{m\cdot S}-457_{m\cdot S},458_{\lambda\cdot r\cdot T},459_{\lambda\cdot r\cdot T},\vphantom{L^{L^{L^{L}}}} \\
465_{W},460_{m\cdot r\cdot T}-492_{m\cdot r\cdot T}\vphantom{L^{L^{L^{L}}}}
\end{array}
& 434\phantom{\centerdot} & 493 & 54\% \\ \hline
20 & 381 &
\begin{array}{c}
381_{a},\overline{\emph{382}}_{t},458_{\lambda},459_{\lambda},465_{W},460_{m\cdot S}-481_{m\cdot S},
482_{\lambda\cdot r\cdot T},\vphantom{L^{L^{L^{L}}}} \\
558_{W},483_{m\cdot r\cdot T}-566_{m\cdot r\cdot T}\vphantom{L^{L^{L^{L}}}}
\end{array}
&458\phantom{\centerdot} & 567 & 59\% \\ \hline
21 & 421 & \overline{\emph{422}}_{t},481_{\lambda},482_{\lambda},483_{m\cdot
S},558_{W},640_{X},484_{m\cdot r}-666_{m\cdot r}\vphantom{L^{L^{L^{L}}}} &
481\phantom{\centerdot} & 667 & 76\% \\ \hline
22 & 463 &
\begin{array}{c}
\overline{\emph{463}}_{br},\overline{\emph{464}}_{t},504_{\lambda},505_{\lambda},558_{W},
506_{m\cdot r}-573_{m\cdot r},574_{\lambda\cdot r},\vphantom{L^{L^{L^{L}}}} \\
640_{X},575_{m\cdot r}-712_{m\cdot r}\vphantom{L^{L^{L^{L}}}}
\end{array}
& 504 & 713 & 84\% \\ \hline
23 & 507 &
\begin{array}{c}
\overline{\emph{507}}_{br},\overline{\emph{508}}_{t},528_{b\cdot f\cdot r},529_{m\cdot r}-553_{m\cdot r},558_{g\cdot W},\vphantom{L^{L^{L^{L}}}} \\
573_{\lambda},574_{\lambda},575_{m},576_{m\cdot r}-744_{m\cdot r}\vphantom{L^{L^{L^{L}}}}
\end{array}
& 573 & 745 & 84\% \\ \hline
24 & 553 &
\begin{array}{c}
553_{a},\overline{\emph{554}}_{t},589_{g},598_{\lambda},599_{\lambda},600_{m\cdot r}-673_{m\cdot r},
\vphantom{L^{L^{L^{L}}}}\\
674_{\lambda\cdot r},675_{m\cdot r}-850_{m\cdot r}\vphantom{L^{L^{L^{L}}}}
\end{array}
 & 598 & 851 & 85\% \\ \hline
25 & 601 &
\begin{array}{c}
620_{g\cdot h},624_{b\cdot f\cdot P\cdot r},625_{m\cdot r}-651_{m\cdot r},673_{\lambda},
\vphantom{L^{L^{L^{L}}}} \\
674_{\lambda},675_{m},906_{W},912_{X},938_{W},676_{m\cdot r}-960_{m\cdot r}
\vphantom{L^{L^{L^{L}}}}
\end{array}
& 673 & 961 & 88\% \\ \hline
26 & 651 &
\begin{array}{c} 651_{a},\overline{\emph{652}}_{t},700_{\lambda},701_{\lambda},
702_{m\cdot r}-781_{m\cdot r},
782_{r},  \vphantom{L^{L^{L^{L}}}}    \\
783_{m\cdot r}-984_{m\cdot r}\vphantom{L^{L^{L^{L}}}}
\end{array}
&700 & 985 & 85\% \\ \hline
27 & 703 &
\begin{array}{c}
728_{b\cdot f\cdot r},729_{m\cdot r}-757_{m\cdot r},781_{\lambda},782_{\lambda},783_{m\cdot S},
\vphantom{L^{L^{L^{L}}}} \\
784_{m\cdot r}-1065_{m\cdot r},1066_{r\cdot T\cdot Z}-
1072_{r\cdot T\cdot Z},\vphantom{L^{L^{L^{L}}}} \\
1073_{m\cdot S\cdot Z}-1103_{m\cdot S\cdot Z},
1104_{r\cdot T\cdot Z},1105_{\lambda\cdot r\cdot T\cdot Z},\vphantom{L^{L^{L^{L}}}}\\
1106_{\lambda\cdot r\cdot T\cdot Z}  \vphantom{L^{L^{L^{L}}}}
\end{array}
&781 & 1107 & 88\% \\ \hline
28 & 757 &
\begin{array}{c}
757_{a},\overline{\emph{758}}_{t},810_{\lambda},811_{\lambda},812_{m\cdot r}-1065_{m\cdot r},\vphantom{L^{L^{L^{L}}}} \\
1066
_{r\cdot T}-1072_{r\cdot T},1073_{m\cdot S}-1103_{m\cdot S},\vphantom{L^{L^{L^{L}}}} \\
1104_{r\cdot T\cdot Z}-
1109_{r\cdot T\cdot Z},1110_{m\cdot r\cdot S\cdot Z}-1141_{m\cdot r\cdot S\cdot Z},
\vphantom{L^{L^{L^{L}}}} \\
1142_{r\cdot T\cdot Z}-1145_{r\cdot T\cdot Z},1146_{\lambda\cdot r\cdot T\cdot Z},
\vphantom{L^{L^{L^{L}}}} \\
1147_{m\cdot r\cdot S\cdot Z}-1170_{m\cdot r\cdot S\cdot Z}
\vphantom{L^{L^{L^{L}}}}
\end{array}
&810 & 1171 & 87\% \\ \hline
\end{array}
\medskip $

Key to Table 7.1: $a\rightarrow (\ref{eq2_cyclicPG(2,q)})$,
$b\rightarrow ( \ref{eq2_cyclicAG(2,q)}),$ $c\rightarrow
(\ref{eq2_cyclicRuzsa}),$ $ f\rightarrow
(\ref{eq2_q-1-cancel}),$ $g\rightarrow (\ref{eq2_Baer}),$ $
h\rightarrow (\ref{eq2_Baer2}),$ $j\rightarrow
(\ref{eq2_FuLabNabDecomp}),$ $ k\rightarrow
(\ref{eq2_affine_q-1}),$ $\lambda\rightarrow
\eqref{eq2_Balbuena1}-\eqref{eq2_AaParBalb}$, $m\rightarrow
(\ref{eq2_tetaExten} ),P\rightarrow
(\ref{eq3.3_AfPlnPairsWeitghts}),r\rightarrow
(\ref{eq4_AfPlExten}),S\rightarrow (\ref{eq4_RuzExten})$,
$T\rightarrow (\ref{eq4_RuzExten_d=p-1}),W\rightarrow $
Table~4.1, $X\rightarrow $Table~4.2$,$ $ y\rightarrow $
Table~6.1 with $k\leq 15$, $Z\rightarrow $ Table~6.1 with $
k>15$, $br\rightarrow $Bruck-Ryser Theorem, $t\rightarrow
$Theorem \ref{th5_deficiency1}
\end{center}
\newpage

\textbf{Table 7.1} (continue 2). The currently known parameters
of symmetric configurations $v_{k}$ (cyclic and
non-cyclic)\smallskip

\begin{center}
$
\begin{array}{@{}r|r|c|r|r|r}
\hline
&  &  & E(k) &  &  \\
k & P(k) & P(k)\leq v\leq G(k)-1 & \leq ~ & G(k) & \text{filling} \\ \hline
29 & 813 &
\begin{array}{c}
\overline{\emph{814}}_{t},840_{b\cdot f\cdot r},841_{m\cdot r}-871_{m\cdot
r},897_{\lambda },898_{\lambda },899_{m\cdot S},\vphantom{L^{L^{L^{L}}}} \\
900_{m\cdot r}-1057_{m\cdot r},1071_{\lambda },1072_{\lambda },1073_{m\cdot
S}-1103_{m\cdot S},\vphantom{L^{L^{L^{L}}}} \\
1104_{r\cdot T}-1109_{r\cdot T},1110_{m\cdot S}-1141_{m\cdot S},
\vphantom{L^{L^{L^{L}}}} \\
1142_{r\cdot T}-1146_{r\cdot T},1147_{m\cdot S}-1179_{m\cdot S},
\vphantom{L^{L^{L^{L}}}} \\
1180_{r\cdot T\cdot Z}-1183_{r\cdot T\cdot Z},1184_{m\cdot S}-1219_{m\cdot
S},\vphantom{L^{L^{L^{L}}}} \\
1220_{r\cdot T\cdot Z},1221_{m\cdot S}-1246_{m\cdot S}
\vphantom{L^{L^{L^{L}}}}
\end{array}
& 1071 & 1247 & 85\% \\ \hline
30 & 871 &
\begin{array}{c}
871_{a},\overline{\emph{872}}_{t},928_{\lambda },929_{\lambda },930_{m\cdot
r}-1057_{m\cdot r},1108_{\lambda },\vphantom{L^{L^{L^{L}}}} \\
1109_{\lambda },1110_{m\cdot S}-1141_{m\cdot S},1142_{r\cdot T}-1146_{r\cdot
T},\vphantom{L^{L^{L^{L}}}} \\
1147_{m\cdot S}-1179_{m\cdot S},1180_{r\cdot T}-1183_{r\cdot T},
\vphantom{L^{L^{L^{L}}}} \\
1184_{m\cdot S}-1217_{m\cdot S},1218_{r\cdot T}-1220_{r\cdot T},1262_{W,}
\vphantom{L^{L^{L^{L}}}} \\
1221_{m\cdot S}-1360_{m\cdot S}\vphantom{L^{L^{L^{L}}}}
\end{array}
& 1108 & 1361 & 78\% \\ \hline
31 & 931 &
\begin{array}{c}
\overline{\emph{931}}_{br},\overline{\emph{932}}_{t},960_{b\cdot f\cdot
r},961_{m\cdot r}-1057_{m\cdot r},1145_{\lambda },\vphantom{L^{L^{L^{L}}}}
\\
1146_{\lambda },1147_{m\cdot S}-1179_{m\cdot S},1180_{r\cdot T}-\
1183_{r\cdot T},\vphantom{L^{L^{L^{L}}}} \\
1184_{m\cdot S}-1217_{m\cdot S},1218_{r\cdot T}\ -1220_{r\cdot T},
\vphantom{L^{L^{L^{L}}}} \\
1221_{m\cdot S}-1255_{m\cdot S},1256_{r\cdot T},\ 1257_{r\cdot T},
\vphantom{L^{L^{L^{L}}}} \\
1262_{W},1258_{m\cdot S\cdot r\cdot T}-1494_{m\cdot S\cdot r\cdot T}
\vphantom{L^{L^{L^{L}}}}
\end{array}
& 1145 & 1495 & 79\% \\ \hline
32 & 993 &
\begin{array}{c}
993_{a},\overline{\emph{994}}_{t},1023_{b\cdot f\cdot r},1024_{m\cdot
r}-1057_{m\cdot r},1182_{\lambda },\vphantom{L^{L^{L^{L}}}} \\
1183_{\lambda },1184_{m\cdot S}-1217_{m\cdot S},1218_{r\cdot T}\
-1220_{r\cdot T},\vphantom{L^{L^{L^{L}}}} \\
1221_{m\cdot S}-1255_{m\cdot S},1256_{r\cdot T},\ 1257_{r\cdot T},
\vphantom{L^{L^{L^{L}}}} \\
1258_{m\cdot S\cdot r\cdot T}-1293_{m\cdot S\cdot r\cdot T},1294_{r\cdot T},
\vphantom{L^{L^{L^{L}}}} \\
1533_{W},1295_{m\cdot S\cdot r\cdot T}-1568_{m\cdot S\cdot r\cdot T}
\vphantom{L^{L^{L^{L}}}}
\end{array}
& 1182 & 1569 & 73\% \\ \hline
33 & 1057 &
\begin{array}{c}
1057_{a},\overline{\emph{1058}}_{t},1219_{\lambda },1220_{\lambda
},1221_{m\cdot S}-1255_{m\cdot S},\vphantom{L^{L^{L^{L}}}} \\
1256_{r\cdot T},1257_{r\cdot T},1258_{m\cdot S\cdot r\cdot T}-1293_{m\cdot
S\cdot r\cdot T},\vphantom{L^{L^{L^{L}}}} \\
1294_{r\cdot T},1634_{W},1295_{m\cdot r}-1718_{m\cdot r}
\vphantom{L^{L^{L^{L}}}}
\end{array}
& 1219 & 1719 & 75\% \\ \hline
34 & 1123 &
\begin{array}{c}
\overline{\emph{1123}}_{br},\overline{\emph{1124}}_{t},1256_{\lambda
},1257_{\lambda },1258_{m\cdot S}-1293_{m\cdot S},\vphantom{L^{L^{L^{L}}}}
\\
1294_{r\cdot T},1295_{m\cdot r}-1429_{m\cdot r},1430_{r\cdot T}-1434_{\
r\cdot T},\vphantom{L^{L^{L^{L}}}} \\
1634_{W},1800_{X},1435_{m\cdot S}-1876_{m\cdot S}\vphantom{L^{L^{L^{L}}}}
\end{array}
& 1256 & 1877 & 82\% \\ \hline
\end{array}
\medskip $

Key to Table 7.1: $a\rightarrow (\ref{eq2_cyclicPG(2,q)})$,
$b\rightarrow ( \ref{eq2_cyclicAG(2,q)}),$ $c\rightarrow
(\ref{eq2_cyclicRuzsa}),$ $ f\rightarrow
(\ref{eq2_q-1-cancel}),$ $g\rightarrow (\ref{eq2_Baer}),$ $
h\rightarrow (\ref{eq2_Baer2}),$ $j\rightarrow
(\ref{eq2_FuLabNabDecomp}),$ $ k\rightarrow
(\ref{eq2_affine_q-1}),$ $\lambda\rightarrow
\eqref{eq2_Balbuena1}-\eqref{eq2_AaParBalb}$, $m\rightarrow
(\ref {eq2_tetaExten} ),P\rightarrow
(\ref{eq3.3_AfPlnPairsWeitghts}),r\rightarrow
(\ref{eq4_AfPlExten}),S\rightarrow (\ref{eq4_RuzExten})$,
$T\rightarrow (\ref {eq4_RuzExten_d=p-1}),W\rightarrow $
Table~4.1, $X\rightarrow $Table~4.2$,$ $ y\rightarrow $
Table~6.1 with $k\leq 15$, $Z\rightarrow $ Table~6.1 with $
k>15$, $br\rightarrow $Bruck-Ryser Theorem, $t\rightarrow
$Theorem \ref {th5_deficiency1}
\end{center}

\newpage

\begin{center}
\textbf{Table 7.1} (continue 3). The currently known parameters
of symmetric configurations $v_{k}$ (cyclic and non-cyclic)

$
\begin{array}{@{}r|r|c|r|r|r}
\hline
&  &  & E(k) &  &  \\
k & P(k) & P(k)\leq v\leq G(k)-1 & \leq ~ & G(k) & \text{filling} \\ \hline
35 & 1191 &
\begin{array}{c}
\overline{\emph{1192}}_{t},1293_{\lambda },1294_{\lambda },1295_{m\cdot
S},1296_{m\cdot r}-1407_{m\cdot r},\vphantom{L^{L^{L^{L}}}} \\
1433_{\lambda },1434_{\lambda },1435_{m\cdot S}-1471_{m\cdot S}, \\
1472_{r\cdot T}-1475_{r\cdot T},1800_{P},1476_{m\cdot S}-1974_{m\cdot S}
\vphantom{L^{L^{L^{L}}}}
\end{array}
& 1433 & 1975 & 83\% \\ \hline
36 & 1261 &
\begin{array}{c}
1330_{\lambda },1331_{\lambda },1332_{m\cdot r}-1407_{m\cdot
r},1474_{\lambda },\vphantom{L^{L^{L^{L}}}} \\
1475_{\lambda },1476_{m\cdot S}-1513_{m\cdot S},1514_{r\cdot T}-1516_{r\cdot
T},\vphantom{L^{L^{L^{L}}}} \\
1517_{m\cdot S}-1519_{m\cdot S},2000_{X},1520_{m\cdot r\cdot T}-2010_{m\cdot
r\cdot T}\vphantom{L^{L^{L^{L}}}}
\end{array}
& 1474 & 2011 & 82\% \\ \hline
37 & 1333 &
\begin{array}{c}
\overline{\emph{1334}}_{t},1368_{b\cdot f\cdot r},1369_{m\cdot
r}-1407_{m\cdot r},\vphantom{L^{L^{L^{L}}}} \\
1515_{\lambda },1516_{\lambda },1517_{m\cdot S}-1555_{m\cdot S},
\vphantom{L^{L^{L^{L}}}} \\
1556_{r\cdot T}-1557_{r\cdot T},1558_{m\cdot r}-2198_{m\cdot r}
\vphantom{L^{L^{L^{L}}}}
\end{array}
& 1515 & 2199 & 83\% \\ \hline
38 & 1407 &
\begin{array}{c}
1407_{a},1556_{\lambda },1557_{\lambda },1558_{m\cdot S},1559_{m\cdot S}, \\
1560_{m\cdot r\cdot S\cdot T}-1597_{m\cdot r\cdot S\cdot T},1598_{r\cdot T},
\vphantom{L^{L^{L^{L}}}} \\
1599_{m\cdot r\cdot T}-1761_{m\cdot r\cdot T},1762_{r\cdot T},1763_{m\cdot
r}-2292_{m\cdot r}\vphantom{L^{L^{L^{L}}}}
\end{array}
& 1556 & 2293 & 83\% \\ \hline
39 & 1483 &
\begin{array}{c}
\overline{\emph{1483}}_{br},\overline{\emph{1484}}_{t},1597_{\lambda
},1598_{\lambda },1599_{m\cdot S},\vphantom{L^{L^{L^{L}}}} \\
1600_{m\cdot r\cdot T}-1761_{m\cdot r\cdot T},1762_{\mathbf{r\cdot T}
},1763_{m\cdot r}-2504_{m\cdot r}\vphantom{L^{L^{L^{L}}}}
\end{array}
& 1597 & 2505 & 89\% \\ \hline
40 & 1561 &
\begin{array}{c}
\overline{\emph{1562}}_{t},1638_{\lambda },1639_{\lambda },1640_{m\cdot
r}-1761_{m\cdot r},\vphantom{L^{L^{L^{L}}}} \\
1762_{\mathbf{\ }r\cdot T},1763_{m\cdot r}-1921_{m\cdot r},1922_{r\cdot
T}-1926_{r\cdot T},\vphantom{L^{L^{L^{L}}}} \\
1927_{m}-1931_{m},1932_{m\cdot r}-2564_{m\cdot r}\vphantom{L^{L^{L^{L}}}}
\end{array}
& 1638 & 2565 & 92\% \\ \hline
41 & 1641 &
\begin{array}{c}
\overline{\emph{1642}}_{t},1680_{b\cdot r},1681_{m\cdot r}-1723_{m\cdot
r},1761_{\lambda },1762_{\lambda },\vphantom{L^{L^{L^{L}}}} \\
1763_{m},1764_{m\cdot r}-1893_{m\cdot r},1925_{\lambda },1926_{\lambda },
\vphantom{L^{L^{L^{L}}}} \\
1927_{m\cdot S}-1969_{m\cdot S},1970_{r\cdot T}-1973_{r\cdot
T},1974_{m}-2610_{m}\vphantom{L^{L^{L^{L}}}}
\end{array}
& 1925 & 2611 & 92\% \\ \hline
42 & 1723 &
\begin{array}{c}
1723_{a},\overline{\emph{1724}}_{t},1804_{\lambda },1805_{\lambda
},1806_{m\cdot r}-1893_{m\cdot r},\vphantom{L^{L^{L^{L}}}} \\
1972_{\lambda },1973_{\lambda },1974_{m\cdot S}-2017_{m\cdot S},
\vphantom{L^{L^{L^{L}}}} \\
2018_{r\cdot T}-2020_{r\cdot T},2021_{m}-2794_{m}\vphantom{L^{L^{L^{L}}}}
\end{array}
& 1972 & 2795 & 85\% \\ \hline
43 & 1807 &
\begin{array}{c}
\overline{\emph{1807}}_{br},\overline{\emph{1808}}_{t},1848_{b\cdot
r},1849_{m\cdot r}-1893_{m\cdot r},2019_{\lambda },\vphantom{L^{L^{L^{L}}}}
\\
2020_{\lambda },2021_{m\cdot S}-2065_{m\cdot S},2066_{r\cdot T}-2067_{r\cdot
T},\vphantom{L^{L^{L^{L}}}} \\
2068_{m\cdot r}-2485_{m\cdot r},2486_{r\cdot T}-2490_{r\cdot T},
\vphantom{L^{L^{L^{L}}}} \\
2491_{m\cdot S}-2593_{m\cdot S},2594_{r\cdot T}-2595_{r\cdot T},
\vphantom{L^{L^{L^{L}}}} \\
2596_{m\cdot S}-3014_{m\cdot S}\vphantom{L^{L^{L^{L}}}}
\end{array}
& 2019 & 3015 & 86\% \\ \hline
\end{array}
\medskip$

Key to Table 7.1: $a\rightarrow (\ref{eq2_cyclicPG(2,q)})$,
$b\rightarrow ( \ref{eq2_cyclicAG(2,q)}),$ $c\rightarrow
(\ref{eq2_cyclicRuzsa}),$ $ f\rightarrow
(\ref{eq2_q-1-cancel}),$ $g\rightarrow (\ref{eq2_Baer}),$ $
h\rightarrow (\ref{eq2_Baer2}),$ $j\rightarrow
(\ref{eq2_FuLabNabDecomp}),$ $ k\rightarrow
(\ref{eq2_affine_q-1}),$ $\lambda\rightarrow
\eqref{eq2_Balbuena1}-\eqref{eq2_AaParBalb}$, $m\rightarrow
(\ref {eq2_tetaExten} ),P\rightarrow
(\ref{eq3.3_AfPlnPairsWeitghts}),r\rightarrow
(\ref{eq4_AfPlExten}),S\rightarrow (\ref{eq4_RuzExten})$,
$T\rightarrow (\ref {eq4_RuzExten_d=p-1}),W\rightarrow $
Table~4.1, $X\rightarrow $Table~4.2$,$ $ y\rightarrow $
Table~6.1 with $k\leq 15$, $Z\rightarrow $ Table~6.1 with $
k>15$, $br\rightarrow $Bruck-Ryser Theorem, $t\rightarrow
$Theorem \ref {th5_deficiency1} \newpage

\textbf{Table 7.1} (continue 4). The currently known parameters
of symmetric configurations $v_{k}$ (cyclic and
non-cyclic)\smallskip

$
\begin{array}{@{}r|r|c|r|r|r}
\hline
&  &  & E(k) &  &  \\
k & P(k) & P(k)\leq v\leq G(k)-1 & \leq ~ & G(k) & \text{filling} \\ \hline
44 & 1893 &
\begin{array}{c}
1893_{a},\overline{\emph{1894}}_{t},2066_{\lambda },2067_{\lambda
},2068_{m\cdot S}-2113_{m\cdot S},\ \vphantom{L^{L^{L^{L}}}} \\
2114_{r\cdot T},2115_{m}-2301_{m},2302_{r},2303_{m}-2485_{m},\vphantom{L^{L^{L^{L}}}} \\
2486_{r\cdot T}-2490_{r\cdot T},2491_{m\cdot S}-2539_{m\cdot S},
\vphantom{L^{L^{L^{L}}}} \\
2540_{r\cdot T}-2543_{r\cdot T},2544_{m\cdot S}-2593_{m\cdot S},
\vphantom{L^{L^{L^{L}}}} \\
2594_{r\cdot T}-2595_{r\cdot T},2596_{m\cdot S}-2647_{m\cdot S},
\vphantom{L^{L^{L^{L}}}} \\
2648_{r\cdot T}-2649_{r\cdot T},2650_{m\cdot S}-3192_{m\cdot S}
\vphantom{L^{L^{L^{L}}}}
\end{array}
& 2066 & 3193 & 86\% \\ \hline
45 & 1981 &
\begin{array}{c}
\overline{\emph{1982}}_{t},2113_{\lambda },2114_{\lambda },2115_{m\cdot
S},2116_{m\cdot r}-2301_{m\cdot r},\vphantom{L^{L^{L^{L}}}} \\
2302_{r},2303_{m\cdot r}-2485_{m\cdot r},\ 2486_{r\cdot T}-2490_{r\cdot T},
\vphantom{L^{L^{L^{L}}}} \\
2491_{m\cdot S}-2539_{m\cdot S},2540_{r\cdot T}-2543_{r\cdot T},
\vphantom{L^{L^{L^{L}}}} \\
2544_{m\cdot S}-2593_{m\cdot S},2594_{r\cdot T}-2596_{r\cdot T},
\vphantom{L^{L^{L^{L}}}} \\
2597_{m\cdot S}-2647_{m\cdot S},2648_{r\cdot T}-2649_{r\cdot T},
\vphantom{L^{L^{L^{L}}}} \\
2650_{m\cdot r\cdot S\cdot T}-2701_{m\cdot r\cdot S\cdot T},2702_{r\cdot T},
\vphantom{L^{L^{L^{L}}}} \\
2703_{m\cdot r\cdot S\cdot T}-3374_{m\cdot r\cdot S\cdot T}
\vphantom{L^{L^{L^{L}}}}
\end{array}
& 2113 & 3375 & 90\% \\ \hline
46 & 2071 &
\begin{array}{c}
\overline{\emph{2072}}_{t},2160_{\lambda },2161_{\lambda },2162_{m\cdot
r}-2301_{m\cdot r},2302_{r},\vphantom{L^{L^{L^{L}}}} \\
2303_{m\cdot r}-2485_{m\cdot r},2486_{r\cdot T}-2490_{r\cdot T},
\vphantom{L^{L^{L^{L}}}} \\
2491_{m\cdot S}-2539_{m\cdot S},2540_{r\cdot T}-2543_{r\cdot T},
\vphantom{L^{L^{L^{L}}}} \\
2544_{m\cdot S}-2593_{m\cdot S},2594_{r\cdot T}-2596_{r\cdot T},
\vphantom{L^{L^{L^{L}}}} \\
2597_{m\cdot S}-2647_{m\cdot S},2648_{r\cdot T}-2649_{r\cdot T},
\vphantom{L^{L^{L^{L}}}} \\
2650_{m\cdot r\cdot S\cdot T}-2701_{m\cdot r\cdot S\cdot T},2702_{r\cdot T},
\vphantom{L^{L^{L^{L}}}} \\
2703_{m\cdot r\cdot S\cdot T}-3446_{m\cdot r\cdot S\cdot T}
\vphantom{L^{L^{L^{L}}}}
\end{array}
& 2160 & 3407 & 93\% \\ \hline
47 & 2163 &
\begin{array}{c}
\overline{\emph{2163}}_{br},\overline{\emph{2164}}_{t},2208_{b\cdot
r},2209_{m\cdot r}-2257_{m\cdot r},\vphantom{L^{L^{L^{L}}}} \\
2301_{\lambda },2302_{\lambda },2303_{m},2304_{m\cdot r}-2451_{m\cdot
r},2489_{\lambda },\vphantom{L^{L^{L^{L}}}} \\
2490_{\lambda },2491_{m\cdot S}-2539_{m\cdot S},2540_{r\cdot T}-2543_{r\cdot
T},\vphantom{L^{L^{L^{L}}}} \\
2544_{m\cdot S}-2593_{m\cdot S},2594_{r\cdot T}-2596_{r\cdot T},
\vphantom{L^{L^{L^{L}}}} \\
2597_{m\cdot S}-2647_{m\cdot S},2648_{r\cdot T}-2649_{r\cdot T},
\vphantom{L^{L^{L^{L}}}} \\
2650_{m\cdot r\cdot S}-2701_{m\cdot r\cdot S},2702_{r\cdot T},2703_{m\cdot
r}-3608_{m\cdot r}\vphantom{L^{L^{L^{L}}}}
\end{array}
& 2489 & 3609 & 91\% \\ \hline
\end{array}
\medskip $

Key to Table 7.1: $a\rightarrow (\ref{eq2_cyclicPG(2,q)})$,
$b\rightarrow ( \ref{eq2_cyclicAG(2,q)}),$ $c\rightarrow
(\ref{eq2_cyclicRuzsa}),$ $ f\rightarrow
(\ref{eq2_q-1-cancel}),$ $g\rightarrow (\ref{eq2_Baer}),$ $
h\rightarrow (\ref{eq2_Baer2}),$ $j\rightarrow
(\ref{eq2_FuLabNabDecomp}),$ $ k\rightarrow
(\ref{eq2_affine_q-1}),$ $\lambda \rightarrow
\eqref{eq2_Balbuena1}-\eqref{eq2_AaParBalb}$, $m\rightarrow
(\ref {eq2_tetaExten}),P\rightarrow
(\ref{eq3.3_AfPlnPairsWeitghts}),r\rightarrow (
\ref{eq4_AfPlExten}),S\rightarrow (\ref{eq4_RuzExten})$,
$T\rightarrow (\ref {eq4_RuzExten_d=p-1}),W\rightarrow $
Table~4.1, $X\rightarrow $Table~4.2$,$ $ y\rightarrow $
Table~6.1 with $k\leq 15$, $Z\rightarrow $ Table~6.1 with $
k>15$, $br\rightarrow $Bruck-Ryser Theorem, $t\rightarrow
$Theorem \ref {th5_deficiency1}

\newpage

\textbf{Table 7.1} (continue 5). The currently known parameters
of symmetric configurations $v_{k}$ (cyclic and
non-cyclic)\smallskip

$
\begin{array}{@{}r|r|c|r|r|r}
\hline
&  &  & E(k) &  &  \\
k & P(k) & P(k)\leq v\leq G(k)-1 & \leq ~ & G(k) & \text{filling} \\ \hline
48 & 2257 &
\begin{array}{c}
2257_{a},\overline{\emph{2258}}_{t},2350_{\lambda },2351_{\lambda
},2352_{m\cdot r}-2451_{m\cdot r},\vphantom{L^{L^{L^{L}}}} \\
2542_{\lambda },2543_{\lambda },2544_{m\cdot S}-2593_{m\cdot S},
\vphantom{L^{L^{L^{L}}}} \\
2594_{r\cdot T}-2596_{r\cdot T},2597_{m\cdot S}-2647_{m\cdot S},
\vphantom{L^{L^{L^{L}}}} \\
2648_{r\cdot T}-2649_{r\cdot T},2650_{m\cdot r\cdot S\cdot T}-2701_{m\cdot
r\cdot S\cdot T},\vphantom{L^{L^{L^{L}}}} \\
2702_{r\cdot T},2703_{m\cdot r}-2881_{m\cdot r},\  \\
2882_{r\cdot T}-2890_{r\cdot T},2891_{m\cdot S}-3774_{m\cdot S}
\vphantom{L^{L^{L^{L}}}}
\end{array}
& 2542 & 3775 & 88\% \\ \hline
49 & 2353 &
\begin{array}{c}
2400_{b\cdot r},2401_{m\cdot r}-2451_{m\cdot r},2595_{\lambda
},2596_{\lambda },\vphantom{L^{L^{L^{L}}}} \\
2597_{m\cdot S}-2647_{m\cdot S},2648_{r\cdot T}-2649_{r\cdot T},
\vphantom{L^{L^{L^{L}}}} \\
2650_{m\cdot r\cdot S\cdot T}-2701_{m\cdot r\cdot S\cdot T},2702_{r\cdot T},
\vphantom{L^{L^{L^{L}}}} \\
2703_{m\cdot r}-2863_{m\cdot r},2889_{\lambda },2890_{\lambda },
\vphantom{L^{L^{L^{L}}}} \\
2891_{m\cdot S}-2941_{m\cdot S},2942_{r\cdot T}-2949_{r\cdot T},
\vphantom{L^{L^{L^{L}}}} \\
2950_{m\cdot S}-3916_{m\cdot S}\vphantom{L^{L^{L^{L}}}}
\end{array}
& 2889 & 3917 & 86\% \\ \hline
50 & 2451 &
\begin{array}{c}
2451_{a},\overline{\emph{2452}}_{t},2648_{\lambda },2649_{\lambda
},2650_{m\cdot S}-2701_{m\cdot S},\vphantom{L^{L^{L^{L}}}} \\
\ 2702_{r\cdot T},2703_{m\cdot r}-2863_{m\cdot r},2948_{\lambda
},2949_{\lambda },\vphantom{L^{L^{L^{L}}}} \\
2950_{m\cdot S}-3001_{m\cdot S},3002_{r\cdot T}-3008_{r\cdot T},
\vphantom{L^{L^{L^{L}}}} \\
3009_{m\cdot S}-4188_{m\cdot S}\vphantom{L^{L^{L^{L}}}}
\end{array}
& 2948 & 4189 & 83\% \\ \hline
51 & 2551 &
\begin{array}{c}
2701_{\lambda },2702_{\lambda },2703_{m\cdot S},2704_{m\cdot r}-2863_{m\cdot
r}, \\
3007_{\lambda },3008_{\lambda },3009_{m\cdot S}-3061_{m\cdot S},
\vphantom{L^{L^{L^{L}}}} \\
3062_{r\cdot T}-3067_{r\cdot T},3068_{m\cdot S}-4380_{m\cdot S}
\vphantom{L^{L^{L^{L}}}}
\end{array}
& 3007 & 4381 & 83\% \\ \hline
\end{array}
\medskip  $

Key to Table 7.1: $a\rightarrow (\ref{eq2_cyclicPG(2,q)})$,
$b\rightarrow ( \ref{eq2_cyclicAG(2,q)}),$ $c\rightarrow
(\ref{eq2_cyclicRuzsa}),$ $ f\rightarrow
(\ref{eq2_q-1-cancel}),$ $g\rightarrow (\ref{eq2_Baer}),$ $
h\rightarrow (\ref{eq2_Baer2}),$ $j\rightarrow
(\ref{eq2_FuLabNabDecomp}),$ $ k\rightarrow
(\ref{eq2_affine_q-1}),$ $\lambda \rightarrow
\eqref{eq2_Balbuena1}-\eqref{eq2_AaParBalb}$, $m\rightarrow
(\ref {eq2_tetaExten}),P\rightarrow
(\ref{eq3.3_AfPlnPairsWeitghts}),r\rightarrow (
\ref{eq4_AfPlExten}),S\rightarrow (\ref{eq4_RuzExten})$,
$T\rightarrow (\ref {eq4_RuzExten_d=p-1}),W\rightarrow $
Table~4.1, $X\rightarrow $Table~4.2$,$ $ y\rightarrow $
Table~6.1 with $k\leq 15$, $Z\rightarrow $ Table~6.1 with $
k>15$, $br\rightarrow $Bruck-Ryser Theorem, $t\rightarrow
$Theorem \ref {th5_deficiency1}
\end{center}
\newpage

\begin{center}
\textbf{Table 7.2.} The currently known parameters of symmetric
configurations $v_{k}$ (cyclic and non-cyclic). Constructions
are not remarked.\smallskip

$\renewcommand{\arraystretch}{1.0}
\begin{array}{@{}c|r|c|r|r|r}
\hline
&  &  & E(k) &  &  \\
k & P(k) & P(k)\leq v\leq G(k)-1 & \leq ~ & G(k) & \text{filling} \\ \hline
3\centerdot & 7 & 7\vphantom{L^{L^{L^{L}}}} & 7\centerdot & 7 & 100\% \\
\hline
4\centerdot & 13 & 13\vphantom{L^{L^{L^{L}}}} & 13\centerdot & 13 & 100\% \\
\hline
5\centerdot & 21 & 21,\overline{\emph{22}}\vphantom{L^{L^{L^{L}}}} &
23\centerdot & 23 & 100\% \\ \hline
6\centerdot & 31 & 31,\overline{\emph{32}},\overline{\emph{33}},34
\vphantom{L^{L^{L^{L}}}} & 34\centerdot & 35 & 100\% \\ \hline
7 & 43 & \overline{\emph{43}},\overline{\emph{44}},45,48,49,50
\vphantom{L^{L^{L^{L}}}} & 48\phantom{\centerdot} & 51 & 75\% \\ \hline
8 & 57 & 57,\overline{\emph{58}},63-68\vphantom{L^{L^{L^{L}}}} & 63
\phantom{\centerdot} & 69 & 67\% \\ \hline
9 & 73 & 73,\overline{\emph{74}},78,80-88\vphantom{L^{L^{L^{L}}}} & 80
\phantom{\centerdot} & 89 & 75\% \\ \hline
10 & 91 & 91,\overline{\emph{92}},98,107-110\vphantom{L^{L^{L^{L}}}} & 107
\phantom{\centerdot} & 111 & 35\% \\ \hline
11 & 111 & \overline{\emph{111}},\overline{\emph{112}},120-133,135-144
\vphantom{L^{L^{L^{L}}}} & 135\phantom{\centerdot} & 145 & 76\% \\ \hline
12 & 133 & 133,\overline{\emph{134}},135,154-170\vphantom{L^{L^{L^{L}}}} &
154\phantom{\centerdot} & 171 & 52\% \\ \hline
13 & 157 & \overline{\emph{158}},168-183,189,193-212\vphantom{L^{L^{L^{L}}}}
& 193\phantom{\centerdot} & 213 & 68\% \\ \hline
14 & 183 & 183,\overline{\emph{184}},210,222-254\vphantom{L^{L^{L^{L}}}} &
222\phantom{\centerdot} & 255 & 50\% \\ \hline
15 & 211 & \overline{\emph{211}},\overline{\emph{212}},231,238-302
\vphantom{L^{L^{L^{L}}}} & 238\phantom{\centerdot} & 303 & 73\% \\ \hline
16 & 241 & 252,255-354\vphantom{L^{L^{L^{L}}}} & 255
\phantom{\centerdot} & 355 & 89\% \\ \hline
17 & 273 & 273,\overline{\emph{274}},288-307,321-398\vphantom{L^{L^{L^{L}}}}
& \mathbf{321}\phantom{\centerdot} & 399 & 79\% \\ \hline
18 & 307 & 307,340-381,389,391,395-398,401,403-432\vphantom{L^{L^{L^{L}}}} &
\mathbf{403}\phantom{\centerdot} & 433 & 63\% \\ \hline
19 & 343 & \overline{\emph{344}},360-381,434-492\vphantom{L^{L^{L^{L}}}}
&434\phantom{\centerdot} & 493 & 54\% \\ \hline
20 & 381 & 381,\overline{\emph{382}},458-566\vphantom{L^{L^{L^{L}}}} &
458\phantom{\centerdot} & 567 & 59\% \\ \hline
21 & 421 & \overline{\emph{422}},481-666\vphantom{L^{L^{L^{L}}}} & 481
\phantom{\centerdot} & 667 & 76\% \\ \hline
22 & 463 & \overline{\emph{463}},\overline{\emph{464}},504-712
\vphantom{L^{L^{L^{L}}}} &504 & 713 & 84\% \\ \hline
23 & 507 & \overline{\emph{507}},\overline{\emph{508}},528-553,558,573-744
\vphantom{L^{L^{L^{L}}}} & 573 & 745 & 84\% \\ \hline
24 & 553 & 553,\overline{\emph{554}},589,598-850\vphantom{L^{L^{L^{L}}}} &
598 & 851 & 85\% \\ \hline
25 & 601 & 620,624-651,673-960\vphantom{L^{L^{L^{L}}}} & 673 & 961 & 88\% \\
\hline
26 & 651 & 651,\overline{\emph{652}},700-984\vphantom{L^{L^{L^{L}}}} &
700& 985 & 85\% \\ \hline
27 & 703 & 728-757,781-1106\vphantom{L^{L^{L^{L}}}} & 781 & 1107 &
88\% \\ \hline
28 & 757 & 757,\overline{\emph{758}},810-1170\vphantom{L^{L^{L^{L}}}} &
810 & 1171 & 87\% \\ \hline
29 & 813 & \overline{\emph{814}},840-871,897-1057,1071-1246
\vphantom{L^{L^{L^{L}}}} &1071 & 1247 & 85\% \\ \hline
30 & 871 & 871,\overline{\emph{872}},928-1057,1108-1360
\vphantom{L^{L^{L^{L}}}} & 1108 & 1361 & 78\% \\ \hline
31 & 931 & \overline{\emph{931}},\overline{\emph{932}},960-1057,1145-1494
\vphantom{L^{L^{L^{L}}}} & 1145 & 1495 & 79\% \\ \hline
32 & 993 & 993,\overline{\emph{994}},1023-1057,1182-1568
\vphantom{L^{L^{L^{L}}}} & 1182 & 1569 & 73\% \\ \hline
\end{array}
$\newpage

\textbf{Table 7.2 }(continue)\textbf{.} The currently known
parameters of symmetric configurations $v_{k}$ (cyclic and
non-cyclic). Constructions are not remarked.\smallskip\\
$\renewcommand{\arraystretch}{0.85}
\begin{array}{@{}c|r|c|r|r|r@{}}
\hline
k & P(k) & P(k)\leq v\leq G(k)-1 & E(k)\leq ~ & G(k) & \text{filling} \\
\hline
33 & 1057 & 1057,\overline{\emph{1058}},1219-1718\vphantom{L^{L^{L^{L}}}} &
1219 & 1719 & 75\% \\ \hline
34 & 1123 & \overline{\emph{1123}},\overline{\emph{1124}},1256-1876
\vphantom{L^{L^{L^{L}}}} & 1256 & 1877 & 82\% \\ \hline
35 & 1191 & \overline{\emph{1192}},1293-1407,1433-1974
\vphantom{L^{L^{L^{L}}}} & 1433 & 1975 & 83\% \\ \hline
36 & 1261 & 1330-1407,1474-2010\vphantom{L^{L^{L^{L}}}} & 1474 &
2011 & 82\% \\ \hline
37 & 1333 & \overline{\emph{1334}},1368-1407,1515-2198
\vphantom{L^{L^{L^{L}}}} & 1515 & 2199 & 83\% \\ \hline
38 & 1407 & 1407,1556-2292\vphantom{L^{L^{L^{L}}}} & 1556 & 2293 &
83\% \\ \hline
39 & 1483 & \overline{\emph{1483}},\overline{\emph{1484}},1597-2504
\vphantom{L^{L^{L^{L}}}} & 1597 & 2505 & 89\% \\ \hline
40 & 1561 & \overline{\emph{1562}},1638-2564\vphantom{L^{L^{L^{L}}}} &
1638 & 2565 & 92\% \\ \hline
41 & 1641 & \overline{\emph{1642}},1680-1723,1761-1893,1925-2610
\vphantom{L^{L^{L^{L}}}} & 1925 & 2611 & 92\% \\ \hline
42 & 1723 & 1723,\overline{\emph{1724}},1804-1893,1972-2794
\vphantom{L^{L^{L^{L}}}} & 1972 & 2795 & 85\% \\ \hline
43 & 1807 & \overline{\emph{1807}},\overline{\emph{1808}},1848-1893,2019
-3014\vphantom{L^{L^{L^{L}}}} & 2019 & 3015 & 86\% \\ \hline
44 & 1893 & 1893,\overline{\emph{1894}},2066-3192\vphantom{L^{L^{L^{L}}}} &
2066 & 3193 & 86\% \\ \hline
45 & 1981 & \overline{\emph{1982}},2113-3374 \vphantom{L^{L^{L^{L}}}} &
2113 & 3375 & 90\% \\ \hline
46 & 2071 & \overline{\emph{2072}},2160-3446\vphantom{L^{L^{L^{L}}}} &
2160 & 3407 & 93\% \\ \hline
47 & 2163 & \overline{\emph{2163}},\overline{\emph{2164}}
,2208-2257,2301-2451,2489-3608 \vphantom{L^{L^{L^{L}}}} & 2489 & 3609 & 91\%
\\ \hline
48 & 2257 & 2257,\overline{\emph{2258}},2350-2451,2542-3774
\vphantom{L^{L^{L^{L}}}} &2542 & 3775 & 88\% \\ \hline
49 & 2353 & 2400-2451,2595-2863,2889-3916\vphantom{L^{L^{L^{L}}}} & 2889
 & 3917 & 86\% \\ \hline
50 & 2451 & 2451,\overline{\emph{2452}},2648-2863,2948-4188
\vphantom{L^{L^{L^{L}}}} & 2948 & 4189 & 83\% \\ \hline
51 & 2551 & 2701-2863,3007-4380\vphantom{L^{L^{L^{L}}}} & 3007 &
4381 & 83\% \\ \hline
52 & 2653 & \overline{\emph{2654}}, 2754-2863,3066-4540
\vphantom{L^{L^{L^{L}}}} & 3066 & 4541 & 84\% \\ \hline
53 & 2757 & \overline{\emph{2758}}, 2808-2863,3125-4694
\vphantom{L^{L^{L^{L}}}} & 3125 & 4695 & 83\% \\ \hline
54 & 2863 & 2863,\overline{\emph{2864}},3184-4746\vphantom{L^{L^{L^{L}}}} &
3184 & 4747 & 83\% \\ \hline
55 & 2971 & \overline{\emph{2971}},\overline{\emph{2972}}, 3243-5196
\vphantom{L^{L^{L^{L}}}} & 3243 & 5197 & 88\% \\ \hline
56 & 3081 & \overline{\emph{3082}}, 3302-5450\vphantom{L^{L^{L^{L}}}} &
3302 & 5451 & 91\% \\ \hline
57 & 3193 & \overline{\emph{3194}}, 3361-5546\vphantom{L^{L^{L^{L}}}} &
3361 & 5547 & 93\% \\ \hline
58 & 3307 & \overline{\emph{3307}},\overline{\emph{3308}},3420-5702
\vphantom{L^{L^{L^{L}}}} & 3420 & 5703 & 95\% \\ \hline
59 & 3423 & \overline{\emph{3424}}, 3480-3541,3597-5822
\vphantom{L^{L^{L^{L}}}} & 3597 & 5823 & 95\% \\ \hline
60 & 3541 & 3541,\overline{\emph{3542}},3658-3783,3838-6038
\vphantom{L^{L^{L^{L}}}} & 3838 & 6039 & 93\% \\ \hline
61 & 3661 & \overline{\emph{3662}}, 3720-3783,3902-6268
\vphantom{L^{L^{L^{L}}}} & 3902 & 6269 & 93\% \\ \hline
62 & 3783 & 3783,3966-6430\vphantom{L^{L^{L^{L}}}} & 3966 & 6431 &
93\% \\ \hline
63 & 3907 & \overline{\emph{3907}},\overline{\emph{3908}},4030-4161,4219-6782
\vphantom{L^{L^{L^{L}}}} & 4219 & 6783 & 94\% \\ \hline
64 & 4033 & 4095-4161,4286-7054\vphantom{L^{L^{L^{L}}}} & 4286 &
7055 & 94\% \\ \hline
\end{array}
$
\end{center}

\begin{thebibliography}{99}
\bibitem{AFLN-graphs} M. Abreu, M. Funk, D. Labbate, and V.
    Napolitano, On (minimal) regular graphs of girth~6,
    \emph{Australas. J. Combin.} \textbf{35}, 119--132 (2006).

\bibitem{AFLN-ConfigGraphs} M. Abreu, M. Funk, D. Labbate, and
    V. Napolitano, Configuration graphs of neighbourhood
    geometries, \emph{Contrib. Discr. Math}. \textbf{3},
    109--122 (2008).

\bibitem{AFLN-Semiplanes} Abreu, M.,  Funk, M., Labbate, D.,
    Napolitano, V.: Deletions, extension, and reductions of elliptic semiplanes,
     Innov. Incidence Geom. \textbf{11}, 139--155 (2010).

\bibitem{AFLN-CyclSchem} Abreu, M.,  Funk, M., Labbate, D.,
    Napolitano, V.: On the ubiquity and utility of cyclic
    schemes, Australas. J. Combin.
 \textbf{55}, 95--120 (2013).

\bibitem{AfDaZ} V. B. Afanassiev, A. A. Davydov, and V. V. Zyablov, Low
density concatenated codes with Reed-Solomon component codes. In: \emph{\
Proc. XI Int. Symp. on Problems of Redundancy in Information and Control
Syst.}, St.-Petersburg, Russia, pp. 47--51 (2007).
http://k36.org/redundancy2007

\bibitem{AfDaZ-InfProc} V. B. Afanassiev, A. A. Davydov, and V.
    V. Zyablov, Low density parity check codes on bipartite
    graphs with Reed-Solomon constituent codes,
    \emph{Information Processes (Electronic Journal)} \textbf{9}, no. 4, 301--331
    (2009).\\
    http://www.jip.ru/2009/301-331-2009.pdf

\bibitem{ArParBalbNetw2011}     Araujo-Pardo, G., Balbuena, C.:
    Constructions of small regular bipartite graphs of girth 6,
    Networks \textbf{57}, 121--127 (2011).

\bibitem{ArParBalbHegDM2010}     Araujo-Pardo, G., Balbuena,
    C.,  H\'{e}ger, T.: Finding small regular graphs of girths 6, 8 and 12
     as subgraphs of cages, Discrete Math. \textbf{310}, 1301--1306
     (2010).

\bibitem{Baker} R. D. Baker, An elliptic semiplane, \emph{J.
    Combin. Theory, ser. A} \textbf{25}, 193--195 (1978).

\bibitem{BalbuenaSIAM2008}   Balbuena,
    C.: Incidence matrices of
    projective planes and of some regular bipartite graphs of
    girth 6 with few vertices, SIAM J. Discrete Math. \textbf{22},
    1351�-1363 (2008).

\bibitem{BaumertGordon} L. D. Baumert and D. M. Gordon, On the existence of
cyclic difference sets with small parameters, In: \emph{High Primes and
Misdemeanours: Lectures in Honour of the 60th Birthday of Hugh Cowie
Williams }, van der Poorten, A., Stein, A., eds., Amer. Math. Soc., Fields
Institute Communications, \textbf{41}, pp. 61--68, Providence (2004).

\bibitem{Boben-v3} M. Boben, Irreducible $(v3)$ configurations and graphs,
\emph{Discrete Math.} \textbf{307}, 331--344, (2007).

\bibitem{Bose} R. C. Bose, An affine analogue of Singer's
    theorem, \emph{J. Ind. Math. Society} \textbf{6}, 1--15
    (1942).

\bibitem{BrasAmorosSemigroup} M. Bras-Amor\'{o}s and K. Stokes,
    The semigroup of combinatorial configurations,
    \emph{Semigroup forum} \textbf{84}, 91-96 (2012).

\bibitem{CarsDinStef-Reduc-n3} H. G. Carstens, T. Dinski, and E. Steffen,
Reduction of symmetric configurations $n_{3}$, \emph{Discrete Applied Math.}
\textbf{99}, 401--411 (2000).

\bibitem{DC2} J. Coykendall and J. Dover, Sets with few intersection numbers
from Singer subgroup orbits, \emph{Europ. J. Combin.} \textbf{22}, 455--464
(2001).

\bibitem{DFGMPConfigArxiv2012} A. A. Davydov, G. Faina, M. Giulietti, S.
Marcugini, and F. Pambianco, On constructions and parameters of symmetric
configurations $v_{k}$. (2012) http://arxiv.org/abs/arXiv:1203.0709v1

\bibitem{DFGMP-submitted} A. A. Davydov, G. Faina, M.
    Giulietti, S. Marcugini, and F. Pambianco, On constructions
    and parameters of symmetric configurations $v_{k}$, submitted.

\bibitem{DGMP-ACCT2008} A. A. Davydov, M. Giulietti, S. Marcugini, and F.
Pambianco, Symmetric configurations for bipartite-graph codes. In: \emph{\
Proc. XI Int. Workshop Algebraic Comb. Coding Theory, ACCT2008,} Pamporovo,
Bulgaria, pp. 63--69 (2008). http://www.moi.math.bas.bg/acct2008/b11.pdf

\bibitem{DGMP-Petersb2009} A. A. Davydov, M. Giulietti, S. Marcugini, and F.
Pambianco, On the spectrum of possible parameters of symmetric
configurations. In: \emph{Proc. XII Int. Symp. on Problems of Redundancy in
Information and Control Systems}, Saint-Petersburg, Russia, pp. 59--64
(2009). http://k36.org/redundancy2009

\bibitem{DGMP-GraphCodes} A. A. Davydov, M. Giulietti, S.
    Marcugini, and F. Pambianco, Some combinatorial aspects of
    constructing bipartite-graph codes, \emph{Graphs
    Combinatorics }\textbf{29}, 187--212 (2013).

\bibitem{DomingoAmorosPeertoPeer} J. Domingo-Ferrer, M. Bras-Amor\'{o}s, Q.
Wu, and J. Manj\'{o}n, User-private information retrieval based on a
peer-to-peer community, \emph{Data Knowl. Eng.} \textbf{68, }1237--1252
(2009).

\bibitem{Dimit} A. Dimitromanolakis, \emph{Analysis of the
    Golomb Ruler and the Sidon Set Problems, and Determination
    of Large, Near-Optimal Golomb Rulers}, Depart. Electronic
    Comput. Eng. Techn. University of Crete, 2002.\\
    http://www.cs.toronto.edu/$\sim $apostol/golomb/main.pdf

\bibitem{Draka} K. Drakakis, A review of the available construction methods
for Golomb rulers, \emph{Advanc. Math. Commun.} \textbf{3}, 235--250 (2009).

\bibitem{Funk1993} M. Funk, On configurations of type $n_{k}$ with constant
degree of irreducibility, \emph{J. Combin. Theory, ser. A,} \textbf{65},
173--201 (1993).

\bibitem{Funk2008} M. Funk, Cyclic difference sets of positive
    deficiency, \emph{Bull. Inst. Combin. Appl.} \textbf{53},
    47--56 (2008).

\bibitem{FunkLabNap} M. Funk, D. Labbate, and V. Napolitano,
    Tactical (de-)compositions of symmetric configurations,
    \emph{Discrete Math. }\textbf{309}, 741--747 (2009).

\bibitem{GabidISIT} E. Gabidulin, A. Moinian, and B. Honary, Generalized
construction of quasi-cyclic regular LDPC codes based on permutation
matrices. In: \emph{Proc. Int. Symp. Inf. Theory 2006, ISIT 2006}, Seattle,
USA, pp. 679--683 (2006).

\bibitem{GH} A. G\'{a}cs and T. H\'{e}ger, On geometric
    constructions of $(k,g)$-graphs,$\ $ \emph{Contrib. Discr.
    Math.} \textbf{3}, 63--80 (2008).

\bibitem{GrahSloan} R. L. Graham and N. J. A. Sloane, On additive bases and
Harmonious Graphs, \emph{Siam J. Algeb. Discrete Methods} \textbf{1},
382--404 (1980).

\bibitem{Gropp-nk} H. Gropp, On the existence and non-existence of
configurations $n_{k}$, \emph{J. Combin. Inform. System Sci.} \textbf{15},
34--48 (1990).

\bibitem{Gropp-Chemic} H. Gropp, Configurations, regular graphs and chemical
compounds, \emph{J. Mathematical Chemistry} \textbf{11}, 145-153 (1992).

\bibitem{Gropp-nonsim} H. Gropp, Non-symmetric configurations with
deficiencies 1 and 2. In: Barlotti, A., Bichara, A., Ceccherini, P.V.,
Tallini, G. (eds.) \emph{Combinatorics '90: Recent Trends and Applications.
Ann. Discrete Math. }vol. 52, pp. 227--239. Elsevier, Amsterdam (1992).

\bibitem{Gropp-ConfGraph} H. Gropp, Configurations and graps --
    II, \emph{Discrete Math.} \textbf{164}, 155--163 (1997).

\bibitem{Gropp-ConfGeomCombin} H. Gropp, Configurations between geometry and
combinatorics, \emph{Discr. Appl. Math. }\textbf{138}, 79--88 (2004).

\bibitem{Gropp-Handb} H. Gropp, Configurations. In: Colbourn, C.J., Dinitz,
J. (eds.) \emph{The CRC Handbook of Combinatorial Designs,} 2-nd edition,
Chapter VI.7, pp. 353--355. CRC Press, Boca Raton (2007).

\bibitem{Grunbaum} B. Gr\"{u}nbaum, \emph{Configurations of points and line}
. Gradute studies in mathematics vol. 103, American Mathematical Society,
Providence (2009).

\bibitem{OstergModulSidon} H. Haanp\"{a}\"{a}, A. Huima, P.
    R. J. \"{O}sterg\aa rd, Sets in $Z_{n}$ with distinct sums
    of pairs, \emph{Discrete Appl. Math.} \textbf{138}, 99--106
    (2004).

\bibitem{Hirs} J. W. P. Hirschfeld, \emph{Projective Geometries
    over Finite Fields}. Second edition, Oxford University
    Press, Oxford (1998).

\bibitem{CyclDecomp} Q. Huang, Q. Diao, S. Lin, Circulant decomposition:
Cyclic, quasi-cyclic and LDPC codes, In: \emph{Proc.} \emph{2010
International Symposium on Information Theory and its Applications (ISITA)},
pp.383--388.

\bibitem{CyclDecomp2} Q. Huang, Q. Diao, S. Lin, K.
    Abdel-Ghaffar, Cyclic and quasi-cyclic LDPC Codes: new
    developments, In: \emph{Proc. Information Theory and
    Applications Workshop (ITA)}, 2011, pp. 1--10.

\bibitem{KaskiOst} P. Kaski, P. R. J. \"{O}sterg\aa rd, There exists no
symmetric configuration with 33 points and line size 6, \emph{Australasian
Journal of Combin.} \textbf{38}, 273--277 (2007).

\bibitem{KaskOst112-11} P. Kaski and P. R. J. \"{O}sterg\aa rd,
    There are exactly five biplanes with k = 11, \emph{J.
    Combin. Des. }\textbf{16}, 117-127 (2007).

\bibitem{Krcadinac} V. Kr\v{c}adinac, \emph{Construction and
    classification of finite structures by computer.} PhD
    thesis, University of Zagreb, in Croatian. (2004).

\bibitem{QCEncoder} Z.-W. Li, L. Chen, L. Zeng, S. Lin, and W.
    H. Fong, Efficient encoding of quasi-cyclic low-density
    parity-check codes, \emph{IEEE Trans. Commun.
    }\textbf{54}, 71--81 (2006).

\bibitem{LingDTSAP} A. C. H. Ling, Difference triangle sets
    from
    affine planes, \emph{IEEE Trans. Inform. Theory}
    \textbf{48}, 2399--2401 (2002).

\bibitem{Lipman} M. J. Lipman, The existence of small tactical
    configurations. In: \emph{Graphs and Combinatorics,
    Springer Lecture Notes in Mathematics} \textbf{406},
    Springer-Verlag, pp. 319--324 (1974).

\bibitem{Longyear1975} J. Q. Longyear, Tactical constructions,
    \emph{J. Combin. Theory, ser. A}, \textbf{19}, 200--207
    (1975).

\bibitem{Martinetti} V. Martinetti, Sulle configurazioni piane $\mu _{3}$.
\emph{Annali di matematica pura ed applicata }(2) \textbf{15} 1--26
(1887-88).

\bibitem{MePaWolk} N. S. Mendelsohn, R. Padmanabhan, and B.
    Wolk, Planar projective configurati-\\ons I, \emph{Note di
    Matematica }\textbf{7}, 91--112 (1987).\\
    http://siba2.unile.it/ese/issues/1/13/Notematv7n1p91.pdf

\bibitem{BibliogrSidon} K. O'Bryant, A complete annotated
    bibliography of work related to Sidon sequences,
    \emph{Electronic J. Combin. }\textbf{11}, $\#39$ (2004).

\bibitem{Pepe} V. Pepe, LDPC codes from the Hermitian curve, \emph{Des.
Codes Crypt.} \textbf{42}, 303--315 (2007).

\bibitem{Ruzsa} I. Z. Ruzsa, Solving a linear equation in a set of integers
I, \emph{Acta Arithmetica }\textbf{65}, 259--282 (1993).

\bibitem{Shearer-Handb} J. Shearer, Difference triangle sets.
    In: Colbourn, C.J., Dinitz, J. (eds.) \emph{The CRC
    Handbook of Combinatorial Designs,} 2-nd edition, Chapter
    VI.19, pp. 436--440. CRC Press, Boca Raton (2007).

\bibitem{ShearerWebShortest} J. Shearer, Table of lengths of
    shortest known Golomb rulers,\\
    http://www.research.ibm.com/peo\-ple/s/shearer/grtab.html

\bibitem{ShearerWebModulGR} J. Shearer, Modular Golomb
rulers,\\
    http://www.research.ibm.com/people/s/shearer/ mgrule.html\\
http://www.research.ibm.com/people/s/shearer/grule.html

\bibitem{ShearerReport} J. B. Shearer, \emph{Difference
    triangle sets constructions,} IBM Research report,
    RC24623(W0808-045) (2008).

\bibitem{Singer} J. Singer, A theorem in finite projective
    geometry and some applications to number theory,
    \emph{Trans. American Math. Soc. }\textbf{43}, 377--385
    (1938).

\bibitem{StokesPeertoPeer} K. Stokes and M. Bras-Amor\'{o}s, Optimal
configurations for peer-to-peer user-private information retrieval.\emph{\
Comput Math Applications} \textbf{59} (2010), 1568--1577.

\bibitem{StokesBrasAmoros2013} K. Stokes and M.
    Bras-Amor\'{o}s, Linear, non-homogeneous, symmetric
    patterns and prime power generators in numerical semigroups
    associated to combinatorial configurations, \emph{Semigroup
    Forum} \textbf{88} (2014), 11--20.

\bibitem{Swanson} C. N. Swanson, Planar cyclic difference packings, \emph{J.
Combin. Des.} \textbf{8}, 426--434 (2000).
\end{thebibliography}
\end{document}